\documentclass[11pt]{article}
\usepackage[a4paper]{geometry}
\usepackage[T1]{fontenc}
\usepackage{lmodern}
\usepackage[utf8]{inputenc}
\usepackage[datamodel=software]{biblatex}
\usepackage{software-biblatex}
\usepackage{amsmath,amssymb,amsfonts,amsthm}
\usepackage{hyperref,enumerate}
\usepackage{tikz-cd}
\usepackage{mathtools}

\newtheorem{theorem}{Theorem}
\newtheorem{proposition}[theorem]{Proposition}
\theoremstyle{definition}
\newtheorem{corollary}[theorem]{Corollary}

\newtheorem{definition}[theorem]{Definition}

\newtheorem{lemma}[theorem]{Lemma}
\newtheorem{example}[theorem]{Example}
\newtheorem{remark}[theorem]{Remark}
\newtheoremstyle {break}{\topsep}{\topsep}{\itshape}{}{\bfseries}{.}{\newline}{}
\theoremstyle {break}
\newtheorem {algorithm}[theorem]{Algorithm}
\newcommand {\inoutput}[2]{{\bf Input:} #1\\
   {\bf Output:} #2}

\numberwithin{equation}{section}
\numberwithin{theorem}{section}

\renewcommand{\epsilon}{\varepsilon}
\renewcommand{\theta}{\vartheta}
\renewcommand{\geq}{\geqslant}
\renewcommand{\leq}{\leqslant}

\newcommand{\point}[1]{(#1)}

\newcommand{\cconj}[1]{\overline {#1}}  
\newcommand{\mmod}[2]{{#1}_{\operatorname {mod} {#2}}} 
\newcommand{\sign}{\operatorname {sign}}
\newcommand{\N}{\mathbf {Z}^{>0}}
\newcommand{\Ng}{\mathfrak{n}}
\newcommand{\Z}{\mathbf {Z}}
\newcommand{\Q}{\mathbf {Q}}
\newcommand{\R}{\mathbf {R}}
\newcommand{\CC}{\mathbf {C}}
\newcommand{\Hc}{\mathbf {H}}
\newcommand{\Cg}{\mathfrak {C}}
\newcommand{\Dg}{\mathfrak {D}}
\newcommand{\NNg}{\mathfrak {N}}
\newcommand{\Pg}{\mathfrak {P}}

\newcommand{\ag}{\mathfrak {a}}
\newcommand{\bg}{\mathfrak {b}}
\newcommand{\cg}{\mathfrak {c}}
\newcommand{\dg}{\mathfrak {d}}
\newcommand{\fg}{\mathfrak {f}}
\newcommand{\mg}{\mathfrak {m}}
\newcommand{\pg}{\mathfrak {p}}

\newcommand{\I}{\mathcal {I}}
\newcommand{\Ac}{\mathcal {A}}
\newcommand{\Bc}{\mathcal {B}}
\newcommand{\Cc}{\mathcal {C}}
\newcommand{\Fc}{\mathcal {F}}
\newcommand{\Pc}{\mathcal {P}}
\newcommand{\Tc}{\mathcal {T}}
\newcommand{\Oc}{\mathcal {O}}
\newcommand{\Qc}{\mathcal {Q}}
\newcommand{\Sc}{\mathcal {S}}
\newcommand{\OK}{\Oc_K}
\newcommand{\OKzero}{\Oc_{K_0}}
\newcommand{\OKr}{\Oc_\Kr}

\newcommand{\Kr}{{K^r}}
\newcommand{\Krzero}{K^r_0}
\newcommand{\omegar}{\omega_r}
\newcommand{\Phir}{{\Phi^r}}
\newcommand{\phir}{{\varphi^r}}
\newcommand{\diff}{\mathcal {D}}
\newcommand{\diffgen}{\lambda}
\newcommand{\rad}{\operatorname {rad}}
\newcommand{\id}{\mathrm {id}}
\newcommand{\Tr}{\operatorname {Tr}}
\newcommand{\Gal}{\operatorname {Gal}}
\newcommand{\Norm}{\operatorname {N}}
\newcommand{\Cl}{\operatorname {Cl}}
\newcommand{\Nphi}{\Norm_\Phi}
\newcommand{\Nphir}{\Norm_{\Phir}}
\newcommand{\NphirO}{\Norm_{\Phir, \Oc}}
\newcommand{\Mat}{\operatorname {Mat}}
\newcommand{\Gl}{\operatorname {GL}}
\newcommand{\Sl}{\operatorname {SL}}
\newcommand{\GSp}{\operatorname {GSp}_{2g}}
\newcommand{\Sp}{\operatorname {Sp}_{2g}}
\newcommand{\Spfour}{\operatorname {Sp}_4}
\newcommand{\GSpplus}{\GSp^+}
\newcommand{\T}[1]{#1^T}
\newcommand{\lquo}[2]
{\leavevmode\kern-.1em\lower.25ex\hbox{$#2$}
\kern-.1em\backslash\kern-.1em\raise.2ex\hbox{$#1$}}
\newcommand{\modstar}{\operatorname {mod^\times}}

\newcommand{\Mg}{\mathbb{M}}

\newcommand{\Vg}{\mathbb{V}}
\newcommand{\Jg}{\mathbb{J}}
\newcommand{\Mo}{M}
\newcommand{\changewithcolumns}[3]{[#3]^{#1}_{#2}}

\newcommand{\otogleft}[1]{\left[ #1 \right]}

\newcommand{\HOPhiN}{H_{\Oc, \Phi}(N)}
\newcommand{\HOPhione}{H_{\Oc, \Phi}(1)}
\newcommand{\HOKPhione}{H_{\OK, \Phi}(1)}

\title{Schertz style class invariants for higher degree CM fields}
\author{Andreas Enge$^1$ and Marco Streng$^2$}
\date{May 9, 2025}

\begin{document}
\maketitle
\footnotetext[1]{ INRIA, Université de Bordeaux, CNRS, CANARI,
33400 Talence, France \newline
\url{https://www.math.u-bordeaux.fr/\textasciitilde aenge/}\newline
\url{andreas.enge@inria.fr}}
\footnotetext[2]{Universiteit Leiden, Postbus 9512, 2300 RA Leiden,
The Netherlands \newline
\url{http://www.math.leidenuniv.nl/\textasciitilde streng/}\newline
\url{streng@math.leidenuniv.nl}

\textbf {Acknowledgements.}
We thank Damien Robert for useful discussions.
This research was partially funded by ERC Starting Grant ANTICS 278537,
and by the Netherlands Organization for Scientific Research (NWO) Vernieuwingsimpuls.
}

\abstract{
Special values of Siegel modular functions for $\Sp (\Z)$ generate
class fields of CM fields. They also
yield abelian varieties with a
known endomorphism ring.
Smaller alternative values of modular functions that lie
in the same class fields (class invariants) thus help
to speed up the computation of those mathematical objects.

We show that modular functions for the subgroup
$\Gamma^0 (N)\subseteq \Sp(\Z)$ yield
class invariants under some splitting conditions on~$N$,
generalising results due to Schertz from
classical modular functions to Siegel modular functions.
We show how to obtain all Galois conjugates of a class invariant by
evaluating the same modular function in CM period matrices derived from
an \emph{$N$-system}. Such a system consists of quadratic polynomials with
coefficients in the real-quadratic subfield satisfying certain congruence
conditions modulo~$N$.
We also examine conditions under which the minimal polynomial of a class
invariant is real.

Examples show that we may obtain class invariants that are much smaller
than in previous constructions.


}

\section {Introduction}

Starting from a \emph {CM field} of degree $2g$, that is, an
imaginary-quadratic extension $K$ of a totally real number field~$K_0$ of
degree~$g$, the theory of \emph{complex multiplication (CM)}
characterises principally polarised abelian
varieties of dimension~$g$ with endomorphism ring an order~$\Oc$ of~$K$.
Invariants of such varieties are algebraic and lie in the
\emph {Shimura class field} of another CM
field, known as the \emph {reflex field}~$\Kr$.
This class field
is contained in the Hilbert class field of~$\Kr$, its maximal unramified
abelian extension. So computing these algebraic invariants
has the two-fold application of constructing abelian varieties with
properties known in advance, and of explicitly constructing class fields
as a first step towards the Hilbert class field of CM fields.

Concretely, the algebraic numbers are obtained as values of a
\emph {Siegel modular function}~$f$ in a CM period matrix $\tau$ in the
\emph {Siegel half space}~$\Hc_g$, a subspace of the symmetric
$g \times g$-matrices with complex coefficients.
We compute the minimal polynomials of these values,
called \emph {class polynomials}.

In the case of $g=1$, corresponding to elliptic curves, the field of Siegel
modular functions of level 1 is generated by the $j$-function. In the case
of $g=2$, corresponding to abelian surfaces and the main focus of this
article, it is generated by the absolute Igusa invariants $i_1$, $i_2$,
$i_3$ of~\cite {igusa}.
It is well-known that the values $f (\tau)$ of these functions~$f$
in a CM period matrix $\tau$ generate the Shimura class field, and algorithms
as well as implementations yielding the associated class polynomials are
available for $g=1$~\cite{Enge09,EnSu10,Sutherland12,cm} and
$g=2$~\cite{Spallek94,Wamelen99,Dupont06,EnTh14,Streng-runtime,recip-3.4.1,cmh-1.1.1}.

However, the algebraic numbers thus obtained have a rather large height,
that is, lead to class polynomials with large coefficients, which require
a proportionally large precision and (over-)proportionally much running
time for their construction.
An approach pursued with success for $g=1$ is to consider other functions~$f$
which are
modular for congruence subgroups $\Gamma^0 (N)$
of some integer level~$N$. The value $f (\tau)$ then lies in an abelian
class field that is generally larger than the Shimura class field (since
it is related not to the Hilbert class field of~$\Kr$, but to its ray class
field of conductor~$N$). However, under certain conditions, $f (\tau)$ lies
in the Shimura class field; we then call it a \emph {class invariant}.
Compared to the $j$- or Igusa invariants, these class invariants have a
height that is generally smaller by an asymptotically constant factor,
which can be as big as~$72$ \cite [Table~7.1]{EnMo14}, and which
significantly increases the range of feasible fields for CM constructions.

The main tool for proving class invariants when $g=1$ in articles such as
\cite{Gee99,GeSt98,Schertz02,EnSc04,EnMo14} is an explicit version of
Shimura's reciprocity law, which expresses the action of the absolute Galois
group of~$\Q$ on $f (\tau)$ via matrix actions on the function~$f$ and the
argument~$\tau$.
While mathematically satisfying, these approaches pose difficulties from
an algorithmic and implementational point of view. They require a good
understanding of Shimura reciprocity, which in general needs to be applied
twice: First one shows that the action of the Galois group of the ray
class field, where $f (\tau)$ lies \textit {a priori}, over the Shimura
class field is trivial on $f (\tau)$, establishing that $f (\tau)$ is
a class invariant, essentially by proving a new mathematical theorem every
time. In a second step, the action of the Galois group of the Shimura
class field over the field $\Kr$ is made explicit to obtain the Galois
conjugates $f_i (\tau_i)$ of $f (\tau)$ and ultimately the class polynomial,
of which they are the roots. It is not straightforward to distill an
algorithm that given~$K$ returns an integer~$N$, a modular function~$f$
of level~$N$, a period matrix~$\tau$ and a polynomial~$H (X) \in \Kr [X]$
such that $f (\tau)$ is a class invariant and a root of~$H (X)$.

In~\cite {Schertz02} Schertz uses Shimura reciprocity for $g=1$ to derive
a rather general criterion for proving class invariants: Roughly speaking,
when the primes dividing~$N$ split or ramify in~$K$, there is a
quadratic polynomial with coefficients in~$\Z$ and with the same
discriminant as~$K$ that represents~$N$, and for $\tau$ a root of this
polynomial and $f$ a function for $\Gamma^0 (N)$ with rational $q$-expansion
coefficients, the value $f (\tau)$ is a class invariant.
In a sense, Schertz applies Shimura reciprocity once and for all; the
result can then be used, without recourse to Shimura reciprocity, as a
sufficient condition to determine an integer~$N$ and a modular function
$f$ of level~$N$ such that $f (\tau)$ is a class invariant.
His criterion has been applied subsequently to prove the existence of
families of class invariants~\cite {EnSc04,EnSc13,EnMo14}, which were
instrumental in certifying primes of record size with elliptic
curve primality proofs~\cite {EnMo02,FrKlMoWi04,Morain07}.

Another important contribution of Schertz's in~\cite {Schertz02}
is to show that all the Galois conjugates of a class invariant $f (\tau)$,
with $\tau$ derived from a quadratic polynomial as sketched above, can be
obtained as $f (\tau_i)$ for the same~$f$, which makes it easier to write
optimised implementations. Moreover, all the $\tau_i$ are derived from a
system of quadratic polynomials satisfying certain congruence conditions
modulo~$2 N$, a so-called \emph {$N$-system}, which may easily be
obtained algorithmically. Altogether, these advances provide a
comprehensive algorithmic treatment of the problem to compute class
polynomials attached to class invariants, and have enabled push-button
implementations written in C on top of standard multiprecision
floating-point libraries~\cite {cm}.

The present article is a step in the endeavour of generalising this
comprehensive algorithmic approach from~$g=1$ to~$g \geq 2$, that is, from
elliptic curves to abelian surfaces and higher dimensional abelian varieties.
From Shimura's adelic formulation of
reciprocity theory, the second author has obtained an explicit
description of the corresponding matrix actions on modular functions and
period matrices, which has enabled him to obtain examples of class
invariants leading to smaller class polynomials~\cite {Streng12}.
We use this as a starting point to derive analogues of Schertz's results
of~\cite {Schertz02} in the case of~$g \geq 2$.

First of all, we show how to obtain a class invariant from a Siegel
modular function of level~$N$ and a quadratic polynomial that satisfies a
congruence condition modulo~$N$.
More precisely, given a function $f$ modular under $\Gamma^0 (N)$ and with
rational $q$-expansion coefficients, and a quadratic polynomial $A X^2 + B X + C$
with coefficients
in the maximal order~$\OKzero$ of the real-quadratic subfield $K_0$ of~$K$,
with $A$ coprime to~$N$ and $C$ divisible by~$N$,
we construct in Theorem~\ref {th:maintheorem}
a period matrix~$\tau$ such that $f (\tau)$ is
a class invariant.

This result provides a sufficient criterion for obtaining a class invariant
from a modular function of level~$N$, albeit conditional to the existence
of a quadratic polynomial satisfying the congruence conditions modulo~$N$,
which is analysed in~\S\ref {ssec:existenceform}. Again, existence of such
a polynomial depends on the splitting behaviour of the primes dividing~$N$ in~$K$.
More precisely, Theorem~\ref {th:n} shows that a suitable polynomial exists if and
only if the prime ideals dividing $N \OKzero$ either split in~$\OK$,
or ramify and occur in $N \OKzero$ with multiplicity~$1$, and the proof of
the theorem provides an explicit algorithm for obtaining such a polynomial.

In a second step we generalise the notion of an $N$-system, a system of
quadratic polynomials representing a certain class group and satisfying congruence
conditions modulo~$2 N$, in Definition~\ref {def:nsystem}, and show in
Theorem~\ref {th:conjugates} that if an $N$-system exists, it describes
all the Galois conjugates of a class invariant obtained from
Theorem~\ref {th:maintheorem}. Again, there is a constructive proof for
the existence of an $N$-system, which we pursue in
\S\ref {ssec:nsystemspractice}.

Together these results can be summarised as follows (with some of the
notions and notations made precise in later sections):

\begin {theorem}
\label {th:summary}
Let $\Oc$ be an order in a CM field~$K$ of degree~$2g$ that is closed under
complex conjugation and contains~$\OKzero$, and assume that the different
of~$K_0$ is principal.
Let $\Phi$ be a primitive CM type such that there
exists a polarised ideal class for $(\Oc, \Phi)$.
Let $N$ be a positive integer,
coprime to the conductor of~$\Oc$, and let $f$ be a Siegel modular
function of level~$N$ that is the quotient of two modular forms with
rational $q$-expansions and invariant under $\Gamma^0 (N)$.
Assume that every prime ideal of $\OKzero$ dividing $N \OKzero$
is either split in $\OK$, or it is ramified and occurs with multiplicity~$1$
in~$N \OKzero$.

Then there exists a quadratic polynomial $Q = A X^2 + B X + C \in \OKzero [X]$
satisfying $N \mid C$, $\gcd (A, N) = 1$ and $A \gg 0$,
which gives rise to a period matrix $\tau$ with CM by $(\mathcal{O},\Phi)$
(Theorem~\ref {th:n}, Proposition~\ref{prop:symplecticbasisgeneral}).
If $\tau$ is not a pole of~$f$, then $f (\tau)$ is a class invariant
(Theorem~\ref {th:maintheorem}).

Let $\Qc = \{ Q_1, \ldots, Q_h \}$ be an $N$-system with $Q_1 = Q$
as in Definition~\ref {def:nsystem}, which exists according to
Theorem~\ref {th:nsystems} and which can be computed by
Algorithm~\ref {alg:nsystem}.
Then the Galois conjugates of the class invariant~$f (\tau)$ are exactly
the~$f (\tau_i)$, where $\tau_i$ is the period matrix obtained from the
quadratic polynomial~$Q_i \in \Qc$ (Theorem~\ref {th:conjugates}).
\end {theorem}

Unlike some previous work, we do not require that~$\Oc$ be the maximal
order~$\OK$ or that $\OKzero$ have (narrow) class number~$1$.
Instead we make the milder assumptions, satisfied by~$\OK$, that
the order~$\Oc$ is closed under complex conjugation (which is
necessary for the existence of principal polarisations,
cf.~Definition~\ref {def:polarisedclass}), and that it
contains the maximal order~$\OKzero$ of the real subfield (which ensures
that the quadratic polynomials under consideration have coefficients in a Dedekind
ring). The only restrictive condition is that the different of~$K_0$ is
principal, but this is always satisfied for $g \leq 2$ (see the discussion
at the beginning of \S\ref {ssec:ideals}).

The next step is to exhibit families of modular functions of level~$N$
as in Theorem~\ref {th:summary} leading to interesting families of class
invariants. We describe a few constructions in~\S\ref {sec:functions}.
While they yield class invariants with smaller heights than the Igusa
invariants, as can be seen from
the examples in~\S\ref{sec:examples}, the gain is not
as spectacular as in the case of~$g=1$. This may be due to the absence
of a function that could play the role of the Dedekind $\eta$-function
and requires further work.

As an interlude, one may notice that while the $j$- and Igusa invariants
define class fields over~$\Kr$, their class polynomials are actually
defined over the smaller totally real field~$\Krzero$.
We relate this property to $N$-systems in~\S\ref{sec:cconj} and
generalise to $g \geq 2$ criteria under which class invariants for $g=1$
have been shown to yield real class polynomials in~\cite{EnSc04,EnMo14}.
These are also illustrated by examples in~\S\ref{sec:examples}.

Future work of the authors will provide analogous results for Hilbert
modular forms, the grounds for which are already laid in the present
article.

\section{The notion of a class invariant}
\label{sec:cft}

The aim of this section is to collect the well-known definitions and
results on Siegel modular functions, complex multiplication and class
field theory needed to give a precise sense to the following definition.

\begin {definition}
\label {def:classinv}
Let on one hand be $f \in \Fc_N$, the field of Siegel modular functions
of level~$N$ and dimension~$g$ with $q$-expansion coefficients in the
cyclotomic field $\Q (\zeta_N)$;
and let on the other hand $\tau$ be a CM point for $(\Oc, \Phi)$,
where $\Oc$ is an order in a CM field $K$ of degree $2g$ and $\Phi$
a CM type.
Then we call $f (\tau)$ a \emph {CM value} or \emph {special value}
of~$f$.

Generically, the value $f (\tau)$ is then an element of $\HOPhiN$, the
\emph{Shimura ray class field of level~$N$} associated to $\Oc$ and $\Phi$ over the
reflex field $\Kr$ of~$K$ (see Section~\ref{ssec:ray});
if $\Oc$ is the maximal order of~$K$, this class field is a
subfield of the ray class field of conductor~$N$ of~$\Kr$.

If moreover $f (\tau) \in \HOPhione$, then we call it a
\emph {class invariant}, and by its \emph{class polynomial}
we mean its characteristic polynomial
\begin {equation}
\label {eq:classpol}
\prod_{[\ag] \in \Cg_{\Oc, \Phi} (1)}
   \left( X - f (\tau)^{\sigma ([\ag])} \right).
\end {equation}
where $\Cg_{\Oc, \Phi} (1)$ is the CM class group of~\eqref {eq:cmgroup}
and $\sigma$ is the Artin map realising its isomorphism with the\
Galois group of the abelian extension $\HOPhione / \Kr$.
\end {definition}

So CM values of Siegel modular functions of level~$1$ are trivially class
invariants; the interest of the definition stems from the fact that
sometimes, under conditions studied in later chapters of this article,
CM values of higher level functions, which are more plentiful, lie in a
class field of smaller conductor than expected (and then generically, they
generate this class field).

\subsection {Siegel modular functions}
\label {ssec:modfunc}

For a commutative ring $R$, let the \textit {symplectic group} be
\begin {equation}
\label {eq:Sp}
\Sp (R) = \left\{ \Mg \in \Mat_{2 g} (R) : \T {\Mg} \Jg \Mg = \Jg
\right\},
\end {equation}
where
\begin {equation}
\label {eq:J}
\Jg = \begin {pmatrix} 0 & 1 \\ -1 & 0 \end {pmatrix}
= \begin {pmatrix} 0 & \id_g \\ -\id_g & 0 \end {pmatrix}.
\end {equation}
The group $\Gamma = \Sp(\Z)$ acts on the \textit {Siegel space} $\Hc_g$,
the set of symmetric complex $g \times g$ matrices with positive definite
imaginary part, by
\begin {equation}
\label {eq:siegelaction}
\begin {pmatrix} a & b \\ c & d \end {pmatrix} \tau
= (a \tau + b)(c \tau + d)^{-1}
\end {equation}
for $\tau \in \Hc_g$, where $a$, $b$, $c$, $d \in \Mat_g (\Z)$.

For a positive integer~$N$, let $\Gamma(N)$ be the kernel of the surjective
reduction map
$\Sp (\Z) \to \Sp (\Z / N \Z)$.
Assuming $g>1$, a \textit {Siegel modular function} of level $N$ is a meromorphic function
on $\lquo {\Hc_g}{\Gamma(N)}$;
for $g=1$ the additional condition of the function being meromorphic at the cusps is
needed.
It can be written as a quotient
of modular forms that have \textit {$q$-expansions}
\[
\sum_T a_T q_T,
\quad a_T \in \CC,
\quad q_T = e^{2 \pi i \Tr (T \tau) / N},
\]
where $T$ runs over the symmetric matrices in
$\Mat_g \left( \frac {1}{2} \Z \right)$ with integral diagonal entries.
Of special interest in the following is the space~$\Fc_N$ of functions
that can be written as quotients of forms with $a_T \in \Q (\zeta_N)$.

The field $\Fc_1$ of rational Siegel modular functions of level~$1$ is
well-known for $g \leq 2$. If $g=1$, it is the $1$-dimensional rational
function field over~$\Q$ generated by the modular $j$-invariant.
If $g=2$, then it is the rational function field of dimension~$3$ over~$\Q$
generated, for instance, by the first three of the eight displayed quotients
on page~642 of Igusa~\cite{igusa}, or alternatively by the
three invariants $i_1$, $i_2$, $i_3$ of~\cite{Streng-runtime},
which are more efficient to use in computations.

\subsection {CM theory}
\label {ssec:cmtheory}

Throughout the remainder of this article, we use the notations and
definitions of~\cite {Streng12}.
Let $K / \Q$ be a CM field of degree~$2 g$, that is, an imaginary-quadratic
extension of a totally real number field $K_0$ of degree~$g$.
Denote by $\Delta_0$ the discriminant of~$K_0$.

\begin {definition}
\label {def:ideal}
Let $\Oc$ be an order of~$K$.
By a \emph {fractional $\Oc$-ideal} we understand a non-zero,
finitely generated $\Oc$-submodule of~$K$.
We call it a \emph {proper fractional $\Oc$-ideal} if $\Oc$ is its
exact ring of multipliers in~$K$.
\end {definition}

By a \emph{CM type} we understand a vector
$\Phi = (\varphi_1, \ldots, \varphi_g):K\rightarrow \CC^g$ representing
the complex embeddings of $K$ up to complex conjugation.
We denote by $\Kr$ the reflex field, another CM field
associated to $\Phi$ and living in the Galois closure of~$K$ in~$\CC$,
and by $\Krzero$ the totally real subfield of~$\Kr$.
A CM type is \emph{primitive} if it is not induced by a CM type of a subfield
of~$K$.
Either all or no CM types of a given quartic CM field are primitive, so
in the case $g=2$ we may
also use the adjective to characterise the field itself.

\begin {definition}
\label {def:polarisedclass}
Let $\Oc$ be an order of~$K$ containing $\OKzero$,
and let $\Phi$ be a primitive CM type of~$K$.
Let $\bg$ be a proper fractional $\Oc$-ideal in the sense of
Definition~\ref {def:ideal}.
Suppose that there exists a $\xi \in K$ with
$\Phi(\xi)\in\left( i \R^{>0} \right)^g$
such that $\xi\bg$ is the trace dual of the complex
conjugate~$\cconj{\bg}$;
if $\Oc=\OK$, then this last condition is equivalent
to $(\bg \cconj \bg \diff_K)^{-1}=\xi\OK$, where $\diff_K$ is
the different of~$K$.
Then the pair $(\bg, \xi)$ is called a \textit {principally polarised ideal}
for $(\Oc, \Phi)$.

Two such pairs $(\bg, \xi)$ and $(\bg', \xi')$
are called \textit {equivalent}
if there is a $\mu\in K^\times$ such that
$\bg' = \mu \bg$ and $\xi' = (\mu \cconj \mu)^{-1} \xi$.
There are finitely many equivalence classes under this relation;
we call them
\textit {principally polarised ideal classes} of $(\Oc,\Phi)$,
and denote the class of $(\bg, \xi)$ by $[(\bg, \xi)]$
and the set of all classes by $\Tc_{\Oc, \Phi}$.
\end {definition}
Note that the existence of a principally polarised ideal class
implies that $\Oc$ is closed under complex conjugation,
as it is the ring of multipliers of both $\bg$ and $\overline{\bg}$.

In the situation of the definition, the bilinear form
$E_\xi : K \times K \to \Q$, $(x, y) \mapsto \Tr (\xi \cconj x y)$
satisfies
$E_\xi (\bg, \bg) = \Z$.
Identifying $\bg$ via $\Phi$ with a $2 g$-dimensional lattice in $\CC^g$
and extending $E_\xi$ to an $\R$-bilinear form on
$\CC^{g}\times \CC^g$ gives a principal polarisation on the complex torus
$\CC^g / \Phi (\bg)$, which has endomorphism ring $\Oc$.
We say that the resulting principally polarised abelian
variety has \emph{CM by $(\Oc, \Phi)$}.
Since $\Phi$ is a primitive CM type, such polarised
abelian varieties are isomorphic if and only if
the associated principally polarised ideals are equivalent.

For any $n \in \N$, the CM type $\Phi$ induces a $\Q$-linear map
$\Phi:K^n \to \Mat_{g \times n} (\CC)$ given by
\[
\Phi : (x_1, \ldots, x_n) \mapsto
\begin {pmatrix}
\varphi_1 (x_1) & \cdots & \varphi_1 (x_n) \\
\vdots && \vdots \\
\varphi_g (x_1) & \cdots & \phi_g (x_n)
\end {pmatrix}.\]
One may choose a symplectic $\Z$-basis $\Sc = (b_1, \ldots, b_{2 g})$ of
$\bg$, that is, a basis
such that the matrix of $E_\xi$ is~$\Jg$ as in~\eqref {eq:J}.
Let $\Sc_1 = (b_1, \ldots, b_g)$
and $\Sc_2 = (b_{g+1}, \ldots, b_{2g})$.
Then
\begin{equation}\label{eq:tauintermsofbasis}
\tau
= \Phi (\Sc_2)^{-1} \Phi (\Sc_1)
= \big( \Phi (b_{g+1}) | \cdots | \Phi (b_{2 g}) \big)^{-1}
\big( \Phi (b_1) | \cdots | \Phi (b_g) \big)
\end{equation}
is called a \textit {CM point}; it is an element of the Siegel space
$\Hc_g$ defined in \S\ref {ssec:modfunc}.
For a Siegel modular function $f \in \Fc_N$, it therefore makes sense
to consider its \textit {CM value}~$f (\tau)$.

\subsection {Class fields}
\label{ssec:ray}

Dealing with non-maximal orders requires a few precautions, but in a class
field theoretic context, we may avoid the finitely many prime ideals that
pose problems.
The \emph{conductor} of $\Oc$ is the $\Oc$- and $\OK$-ideal
$\fg = \fg_{\Oc} = \{x\in K : x\OK\subseteq \Oc\} \subseteq \Oc$.

The monoid of integral ideals of~$\Oc$ coprime to~$\fg$ is isomorphic
to the monoid of integral ideals of~$\OK$ coprime to~$\fg$ via the map
$\ag \mapsto \ag_K := \ag \OK$
and its inverse $\ag_K \mapsto \ag = \ag_K \cap \Oc$,
see the proof of \cite [Proposition~7.20]{Cox89},
which is formulated for imaginary-quadratic fields,
but carries over immediately to arbitrary number fields.
An integral ideal $\ag$ of~$\Oc$ coprime to~$\fg$, by which we
mean that $\ag + \fg = \Oc$, is invertible,
cf.~\cite[Propositions (12.4) and (12.10)]{Neukirch99}.

Let $F$ be the positive integer such that $\fg \cap \Z = F \Z$.
To simplify, from now on we will assume that all integral or
fractional ideals are coprime to~$F$.
Then additional coprimality conditions in~$\Oc$ can be expressed
in terms of the Dedekind ring $\OK$:
A  non-zero integral ideal~$\ag$ of~$\Oc$ is coprime to~$NF$ if and only if
$\ag = \ag_K \cap \Oc$ for an integral ideal~$\ag_K$ of~$\OK$
such that
$v_{\pg} (\ag_K) = 0$ for all primes $\pg \mid NF$,
and a fractional $\Oc$-ideal~$\cg$ is coprime to~$NF$ if and only if
it can be written as $\cg = \ag \bg^{-1}$ with non-zero integral ideals
$\ag$ and $\bg$ of $\Oc$ that are coprime to~$NF$.

We define coprimality and congruences for elements of $K^\times$
as in~\cite[Definition~4.2]{Streng12}.
\begin{definition}
\label {def:modstar}
For a positive integer $N$ and for $x\in K^\times$,
we say that $x$ is \emph{coprime to $NF$ with respect to $\Oc$}
if one of the following equivalent conditions holds,
where $a \mapsto a'$ denotes reduction modulo~$N F$ in~$\Oc$.
\begin{enumerate}
\item \label{point2} $x = a/b$ for some $a\in \Oc$ and $b\in\Z\setminus\{0\}$ with $a'\in(\Oc/NF\Oc)^\times$
and $b\in 1+NF\Z$;
\item \label{point4} $x\Oc = \ag \bg^{-1}$ for non-zero $\Oc$-ideals $\ag$ and $\bg$ that are
coprime to~$N F$.
\end{enumerate}
We write $x \equiv 1\modstar N\Oc$ to mean that the following condition holds:
\begin{itemize}
\item[(1')] as in \point{\ref{point2}} above, with additionally $a-1\in N\Oc$.
\end{itemize}
\end{definition}

Now let us return to our specific situation of a CM field~$K$
and an order~$\Oc$ of $K$ with $\cconj \Oc = \Oc$
and such that $\Oc$ contains $\OKzero$.
Denote by $\I (NF)$ the group of fractional ideals of $\OKr$
that are coprime to~$NF$.
To a CM type $\Phi$ of $K$ one may associate a \emph {reflex CM type}
$\Phir$ of $\Kr$. Then the reflex type norm is
the multiplicative map $\Kr \to K$ given by
$\Nphir (\alpha) = \prod_{\phir \in \Phir} \phir (\alpha)$.
It extends naturally to a map on ideals, which sends ideals
of $\Oc_{K^r}$ that are coprime to $NF$ to
ideals of $\OK$ that are coprime to~$NF$.
Intersecting with~$\Oc$ leads to ideals of~$\Oc$ coprime to~$N F$,
and we denote the resulting map by $\NphirO$.
Extending multiplicatively, we get a homomorphism $\NphirO$ from
the group
$\I(NF)$ of fractional $\OKr$-ideals coprime to~$NF$
to the group of fractional $\Oc$-ideals coprime to~$NF$.

Let the \emph {CM class group} for $(\Oc, \Phi)$ of level~$N$ be
defined by
\begin {equation}
\label {eq:cmgroup}
\Cg_{\Oc, \Phi} (N) = \I (NF) / S_{\Oc, \Phi} (N),
\end {equation}
where
\begin {equation}
\label {eq:ray}
\begin {split}
S_{\Oc, \Phi} (N) = \{ \ag \in \I (N F)\ :\ \
& \exists \mu \in K^\times \text { with }
\NphirO(\mathfrak{a}) = \mu\mathcal{O},\ \mu\cconj{\mu}\in\Q, \\
& \mu \equiv 1 \modstar N \Oc \}.
\end {split}
\end {equation}
The CM class group of level~$1$, $\Cg_{\Oc, \Phi} (1)$, acts freely on the
set $\Tc_{\Oc, \Phi}$ of principally polarised ideal classes as given in
Definition~\ref {def:polarisedclass}, via
\begin {equation}
\label {eq:cmgroupaction}
[\ag] \cdot [(\bg, \xi)]
= \left[ ( \NphirO (\ag)^{-1} \bg, \Norm (\ag) \xi ) \right],
\end {equation}
where $[\ag]$ denotes the class of $\ag \in \I (F)$.
So the role of the type norm map is to connect the realms of CM class groups
of~$\Kr$ and of principally polarised ideal classes of~$K$.

For a fixed CM point $\tau$ derived from a principally polarised ideal class in
$\Tc_{\Oc, \Phi}$, let $\HOPhiN\subseteq \CC$ be the field generated
over~$\Kr$
by all values $f (\tau)$ for the $f \in \Fc_N$ that are regular at~$\tau$.
Then $\HOPhiN$ is, independently of $\tau$,
the abelian class field of~$\Kr$ with Galois group isomorphic to
$\Cg_{\Oc, \Phi} (N)$, see \cite [Theorem~2.5]{Streng12}
or \cite[Main theorem~3, p.~142]{shimura-taniyama}.
We call $\HOPhiN$ the \emph{Shimura ray class field of level~$N$},
or if $N=1$ simply the \emph{Shimura class field}.
We denote the \emph {Artin map}, which realises the isomorphism
between the class group and the Galois group, by
\begin{equation}\label{eq:artinmap}
\sigma = \sigma_N : \Cg_{\Oc, \Phi} (N) \xrightarrow{\sim}
\Gal (\HOPhiN / \Kr).
\end{equation}
As $S_{\Oc, \Phi} (N)$ contains the principal ray of modulus $N F$, the
field $\HOPhiN$ is a subfield of the ray class field of
modulus~$N F$ of~$\Kr$.
In particular, the field $\HOKPhione$
is a subfield of the Hilbert class field of~$K^r$.

\section {Class invariants from functions for
\texorpdfstring {$\Gamma^0 (N)$}{Gamma0 (N)}}
\label{sec:classinv}

\subsection {Polarised ideal classes and symplectic bases}
\label {ssec:ideals}

To get an explicit handle on ideals and polarised ideal classes,
we would like to mimic the situation for $g = 1$, where ideals are
represented as $z \Z + \Z$ with $z \in K$.
In higher degree CM fields, one would hope for an analogous representation
with $\OKzero$ in the place of~$\Z$. While such a representation need not
exist in general, it does in a more special situation we assume from now on:
Let $\Oc$ be an order in a CM field $K$ such that $\Oc\supseteq \OKzero$
and $\diff_{K_0}=\diffgen\OKzero$ for some $\diffgen\in K_0$,
where $\diff_{K0}$ is the different of~$K_0$.
Such a $\diffgen$ exists, for instance, when $g=1$ (with $\diffgen=1$)
and when $g=2$ (with $\diffgen=\sqrt{\Delta_0}$,
the square root of the discriminant of~$K_0$).

We then get the following classification of principally polarised ideal
classes.
\begin {proposition}
\label {prop:twogenerators}
Let $K$ be a CM field such that $\diff_{K_0} = \diffgen \OKzero$
for some $\lambda \in K_0$.

To $z\in K\setminus K_0$ associate $\bg = z \OKzero + \OKzero$ and
$\xi = \left( (z - \cconj z) \diffgen \right)^{-1}$.
Then $\Oc = \{x\in K : x\bg\subseteq \bg\}$ is an order in $K$ that
is stable under complex conjugation and contains~$\OKzero$,
the set $\Phi = \{\varphi : K\rightarrow \CC \mid \varphi(\xi)\in i\R^{>0}\}$
is a CM type of~$K$, and the pair
$(\bg, \xi)$ is a principally polarised ideal for~$(\Oc, \Phi)$.

Conversely, every principally polarised ideal class for
any CM type $\Phi$ of $K$ and any order
$\Oc$ of $K$ that is stable under complex conjugation and contains
$\OKzero$ has such a representative.
\end {proposition}

\begin {proof}
The results are \cite [Theorems~I.5.8--9]{Streng10},
where they are stated for maximal orders $\Oc$,
but the proof only uses that $\OKzero$ is maximal.
\end {proof}

We may then write down an explicit symplectic basis
and period matrix for such a polarised ideal.

\begin {proposition}
\label {prop:symplecticbasisgeneral}
Assuming that $\diff_{K_0} = \diffgen\OKzero$, let $z$, $\bg$, $\xi$
and $\Phi$ be as in Proposition~\ref{prop:twogenerators}.
Let $\Bc_1 =(b_{1,1}, \ldots, b_{1,g})$ be any $\Z$-basis of $\OKzero$.
Write its trace-dual $\Q$-basis of~$K_0$ as
$ -\diffgen^{-1} \Bc_2
= (-\diffgen^{-1} b_{2,1},\ldots, -\diffgen^{-1} b_{2,g})$.
Then a symplectic basis of $(\bg, \xi)$ is given by
$$\Sc = (z b_{1,1},\ldots, z b_{1,g}, b_{2,1},\ldots, b_{2,g})
= (z \Bc_1 | \Bc_2),$$
and a period matrix by $\tau = \Phi(\Bc_2)^{-1}\Phi(z\Bc_1)$.
\end{proposition}
\begin{proof}
Note first that the trace-dual is a $\Z$-basis of $\diff_{K_0}^{-1}$, so
$\Sc$ is indeed a basis of $\bg = z\OKzero + \OKzero$.

Next note that since $\xi$ is purely imaginary, that is,
$\cconj \xi = - \xi$, we have
$\Tr (\xi \alpha)
= \Tr_{K_0 / \Q} \big( \xi (\alpha - \cconj \alpha) \big)$
for any $\alpha \in K$.
Since for $(u,v) = (z b_{1,i}, z b_{1,j})$
and for $(u,v) = (b_{2,i}, b_{2,j})$ we have $\cconj{u}v\in K_0$,
this implies
$E_{\xi} (u,v) = \Tr (\xi \cconj{u} v) = 0$.
Finally,
\begin{align*}
E_{\xi} (z b_{1,i}, b_{2,j}) &=
\Tr_{K_0 / \Q} \left(
(z - \cconj z)^{-1} \diffgen^{-1}
(\cconj z b_{1, i} b_{2, j} - z b_{1, i} b_{2, j}) \right) \\
&=
\Tr_{K_0/\Q}(- \diffgen^{-1} b_{2,j} b_{1,i}) = \delta_{i j},
\end{align*}
hence the basis is symplectic.
The formula for the period matrix is~\eqref{eq:tauintermsofbasis}.
\end{proof}

\begin {corollary}
\label {cor:symplecticbasisg2}
Let $g = 2$ and $\diffgen = \sqrt {\Delta_0}$.
In the situation of Proposition~\ref {prop:twogenerators},
let $(\varphi_1,\varphi_2) = \Phi$, and to simplify
the notation, write $\alpha_i = \varphi_i (\alpha)$ for any $\alpha\in K$.
A symplectic basis $\Sc$ of $\bg$ with respect to $E_\xi$ and
an associated period matrix~$\tau$ are given as follows:

If $\Delta_0$ is odd, let $\omega = \frac {1 + \diffgen}{2}$;
then $\Sc = (z \omega, z,-1, 1-\omega)$.
If $\Delta_0$ is even, let $\omega = \frac {\diffgen}{2}$;
then $\Sc = (z \omega, z,-1, -\omega)$.
In both cases,
\[
\tau = \frac {1}{- \lambda_1}
\begin {pmatrix}
z_1 \omega_1^2 - z_2 \omega_2^2 & z_1 \omega_1 -  z_2 \omega_2 \\
z_1 \omega_1 - z_2 \omega_2     & z_1 - z_2
\end {pmatrix}
.
\]
\end {corollary}
\begin{proof}
Take $\Bc_1=(\omega, 1)$, $b_{2,1} = -1$, and $b_{2,2} = -\omega$
if $\Delta_0$ is even
or $b_{2,2} = 1-\omega$ if $\Delta_0$ is odd.
It is easy to check that
$(-\diffgen^{-1} b_{2,1}, -\diffgen^{-1} b_{2,2})$
is the trace dual basis of~$\Bc_1$, so the result follows from
Proposition~\ref{prop:symplecticbasisgeneral} using
$\lambda_2 = - \lambda_1$.
\end{proof}

\subsection {Quadratic polynomials and polarised ideal classes}
\label {ssec:quadpols}

We have seen in Proposition~\ref {prop:twogenerators} that a polarised
ideal class can be represented by a pair $(\bg, \xi)$, in which the
fractional ideal~$\bg$ satisfies $\bg = z \OKzero + \OKzero$ for some
$z \in K \backslash K_0$, and $\xi \in K$ is also computed as a function
of~$z$.
So ultimately~$z$ determines the class.
It is conveniently given as the root of an irreducible quadratic polynomial
$A X^2 + B X + C$ with coefficients $A, B, C \in \OKzero$.

\begin {definition}
\label {def:represents}
Let $T \in \Tc_{\Oc, \Phi}$ be a principally polarised ideal class
for $(\Oc, \Phi)$ and let $Q = A X^2 + B X + C \in \OKzero [X]$.
Then we say that $Q$ \emph {represents} $T$ if there is a root~$z$
of~$Q$ such that $T = [(\bg, \xi)]$ with $\bg$ and $\xi$
obtained from~$z$
by Proposition~\ref {prop:twogenerators}.
\end {definition}

Note that if $Q$ represents $T$, then only one root $z$ of the two
roots of $T$ leads to a principally polarised ideal for $(\Oc, \Phi)$:
The other one, $\cconj {z}$, leads to $\cconj {\xi} = - \xi$, which
is purely negative imaginary under the embeddings in $\Phi$ instead of
purely positive imaginary. So when the context is clear, we call this $z$
\emph {the} root of~$Q$.
(Another way to look at this is to notice that $\cconj {z}$ leads by
Proposition~\ref {prop:twogenerators} to the same order $\cconj {\Oc} = \Oc$
and to the CM type $\cconj {\Phi}$, which is equivalent to $\Phi$
by CM theory;
or without recourse to theory, one notices that replacing $z$ by
$\cconj {z}$ and $\Phi$ by $\cconj {\Phi}$ in
Proposition~\ref{prop:symplecticbasisgeneral}
and Corollary~\ref{cor:symplecticbasisg2} leads to the same value
of~$\tau$.)

For $g=1$, one usually assumes $Q$~to be primitive, that is, with coprime
coefficients in~$\Z$ and $A>0$; this choice makes the polynomial unique.
Unless the narrow class number of~$K_0$ is~$1$, we cannot hope to achieve
this in general, so we need to adopt a weaker convention:
We may at least avoid any finite set of primes in the greatest common
divisor; and, more strongly, we will see that each polarised ideal class has a
representative in
which the coefficient~$A$ is not divisible by any of these primes
(Proposition~\ref{prop:Acoprime}).
We are mainly interested in congruences of the
quadratic polynomial modulo rational integers; but since the proofs are
identical, we formulate generalisations modulo ideals $\Ng$ of~$\OKzero$,
which most of the time will be $\Ng = N \OKzero$ for a rational
integer~$N$.

In a first step, let us consider a notion weaker than primitivity, which
can be made to hold for arbitrary $z \in K\setminus K_0$.

\begin{definition}
\label{def:semiprimitive}
Let $\Ng$ be a non-zero ideal of $\OKzero$.
A quadratic polynomial $AX^2+BX+C\in\OKzero[X]$ is \emph{semiprimitive
modulo $\Ng$}
if $A$ is totally positive and furthermore $\gcd(A,B,C,\Ng)=1$.
\end{definition}

\begin {proposition}
\label {prop:quadformexists}
Let $\Ng$ be a non-zero ideal of $\OKzero$.
Every element $z$ of $K \backslash K_0$ is a root of a
quadratic polynomial $A X^2 + B X + C\in\OKzero[X]$
that is semiprimitive modulo~$\Ng$.
The discriminant $B^2-4AC$ of the polynomial is totally negative.
The $\OKzero$-module $\bg = z \OKzero + \OKzero$
is an invertible (hence proper) fractional ideal of the order
\begin {equation}
\label {eq:orderfromform}
\Oc = \dg^{-1} A z + \OKzero,
\end {equation}
where $\dg = \gcd (A, B, C) = A \OKzero + B \OKzero + C \OKzero$.
\end{proposition}

Definition~\ref{def:semiprimitive} and Proposition~\ref{prop:quadformexists}
do not need $\diff_{K_0}$ to be principal, even though the rest of
Section~\ref{sec:classinv} does.

\begin {proof}
As $K=K_0(z)\supseteq K_0$ is a quadratic extension,
there is a non-zero polynomial $AX^2+BX+C\in K_0[X]$
with~$z$ as a root.
By the strong approximation theorem, for instance
\cite [Corollary~1.2.9]{Cohen00},
there is an element $d \in K_0$ such that $v_\pg (d) = - v_\pg (\dg)$
for each prime ideal $\pg$ dividing~$\Ng$,
such that $v_\pg (d) \geq 0$
for all other non-zero prime ideals~$\pg$, and such that
the signs of $d$ under the real embeddings
of~$K_0$ coincide with those of~$A$. Then we may multiply $A$, $B$ and $C$
by $d$ to obtain new coefficients in $\OKzero$ with $\gcd(A,B,C)$ coprime to~$\Ng$
and $A$ totally positive.

The discriminant is totally negative
as $K=K_0(z)\supsetneq K_0$ is totally imaginary-quadratic.

The set $\Oc$ is a subring of $K$ as $A^2z^2 = -ABz-AC\in \dg^2\Oc$.
The $\OKzero$-module $\bg$ is a fractional $\Oc$-ideal as both $Az$ and $Az^2 = -Bz-C$ lie in $\dg\bg$.
Moreover, we have $\bg\cconj{\bg} = z\cconj{z}\OKzero + z\OKzero + \cconj{z}\OKzero + \OKzero
= (C/A)\OKzero + z\OKzero + (B/A)\OKzero + \OKzero = \dg\Oc/A$,
so $\bg$ is an invertible fractional $\Oc$-ideal with inverse $\dg^{-1}A\cconj{\bg}$.
\end{proof}

To go further, we can use the leeway provided by the possibility of
changing the representative of a polarised ideal class. So we need
to consider under which conditions two numbers
$z, z' \in K \backslash K_0$ represent the same polarised ideal class.

\begin {proposition}
\label {prop:equivz}
Assuming again $\diff_{K_0} = \diffgen\OKzero$,
associate to $z, z' \in K \backslash K_0$ the principally polarised
ideals $(\bg, \xi)$ for $(\Oc, \Phi)$ and $(\bg', \xi')$ for $(\Oc', \Phi')$
as in Proposition~\ref {prop:twogenerators}.
Then the following assertions are equivalent:
\begin{enumerate}
\item
   We have $(\Oc', \Phi') = (\Oc, \Phi)$ and $[(\bg',\xi')] = [(\bg,\xi)]$.
\item There is a matrix
\[
\Mo =
\begin {pmatrix}
\alpha & \beta \\ \gamma & \delta
\end {pmatrix}
\in \Sl_2 (\OKzero)\]
such that
$z' = \Mo z := \frac {\alpha z + \beta}{\gamma z + \delta}$.
\end{enumerate}
If the equivalent assertions hold, then we have
$\xi' = (\gamma z + \delta)(\gamma \cconj z + \delta) \xi$.
If furthermore $z$ is a root of $A X^2 + B X + C \in \OKzero [X]$,
then $z'$ is a root of $A' X^2 + B' X + C'\in \OKzero [X]$ with
\begin {equation}
\label {eq:equivform}
\begin {split}
A' & = A \delta^2 - B \gamma \delta + C \gamma^2, \\
B' & = - 2 A \beta \delta + B (1 + 2 \beta \gamma)
       - 2 C \alpha \gamma, \\
C' & =  A \beta^2 - B \alpha \beta + C \alpha^2,
\end {split}
\end {equation}
and
\begin {equation}
\label {eq:equivxi}
\xi' = \frac {A'}{A} \, \xi.
\end {equation}
Finally, we have $\gcd (A', B', C') = \gcd (A, B, C)$, and $A'$ is totally
positive if $A$ is, so that semiprimitivity of the polynomial for~$z$
carries over to the polynomial for~$z'$.
\end {proposition}

\begin {proof}
Given any
\[\Mo = \begin {pmatrix} \alpha & \beta \\ \gamma & \delta \end {pmatrix}
\in \Gl_2 (\OKzero),\]
one easily computes
\begin{equation}\label{eq:sl2k0}
\frac{z-\cconj{z}}{\Mo z-\cconj{\Mo z}}
= (\gamma z+\delta)(\gamma\cconj{z}+\delta) (\det \Mo)^{-1}.
\end{equation}

Assume first that there exists $\Mo\in \Sl_2(\OKzero)$
with
$z' = \Mo z$.
Then one has $\bg = (\alpha z + \beta) \OKzero + (\gamma z + \delta)
\OKzero$
and $\bg' = z' \OKzero + \OKzero = \mu \bg$
with $\mu = (\gamma z + \delta)^{-1}$.
This implies $\Oc' = \Oc$.
By~\eqref{eq:sl2k0} one sees that
$\xi' = \xi (\gamma z + \delta)(\gamma \cconj z + \delta)
= \xi (\mu \cconj \mu)^{-1}$.
So $(\bg, \xi)$ and $(\bg', \xi')$ are indeed equivalent,
and $\xi'$ belongs to the same CM type $\Phi$ as $\xi$
since $\mu \cconj \mu$ is totally positive.

Conversely, if the two pairs are equivalent for the same $(\Oc, \Phi)$,
then $\bg' = \mu \bg$ and $\xi'=(\mu\cconj{\mu})^{-1}\xi$
for some $\mu \in K^\times$, which implies
$z' = \mu (\alpha z + \beta)$ and $1 = \mu (\gamma z + \delta)$
for some $\alpha$, $\beta$, $\gamma$, $\delta \in \OKzero$,
so that $z' = \Mo z$ with
\[\Mo =
\begin {pmatrix}
\alpha & \beta \\ \gamma & \delta
\end {pmatrix}
\in \Gl_2 (\OKzero)\]
as the transformation is invertible.
Now the definition of $\xi$ and $\cconj \xi$ and~\eqref{eq:sl2k0}
yield $1 = \xi / (\xi' \mu\cconj{\mu}) = \det \Mo$, so
$\Mo\in\Sl_2(\OKzero)$.

Now define $A'$, $B'$, and $C'$ by~\eqref{eq:equivform}.
A direct verification shows that
$z' = \Mo z$ is a root of $A'X^2+B'X+C$
(in fact, we have $A'X^2 + B'XZ + C'Z^2 = (AX^2+BXZ+CZ^2) \circ \Mo^{-1}$).
We also have \[\frac{\xi'}{\xi} = (\gamma z + \delta)(\gamma\overline{z}+\delta)
= \gamma^2 \frac{C}{A} + \gamma\delta\frac{-B}{A} + \delta^2 = \frac{A'}{A},\]
which is \eqref{eq:equivxi} and proves that $A'$ is totally positive if $A$ is.

From~\eqref {eq:equivform} one reads off that
$\gcd (A, B, C) \mid \gcd (A', B', C')$;
noticing that the inverse matrix $\Mo^{-1}$ leads to similar formul{\ae}
to express $(A', B', C')$ in terms of $(A, B, C)$ shows the converse.
\end {proof}

Now we have all ingredients to formulate and prove the main result of
Section~\ref{ssec:quadpols}.

\begin {proposition}
\label {prop:Acoprime}
Let $\Ng$ be a non-zero ideal of $\OKzero$, and assume that
$\Oc \supseteq \OKzero$ and $\diff_{K_0} = \diffgen \OKzero$.
Then every principally polarised ideal class for $(\Oc, \Phi)$
is represented by a polynomial $A X^2 + B X + C \in \OKzero [X]$
with $A$ totally positive and $\gcd (A \OKzero, \Ng) = \OKzero$.
\end {proposition}

\begin {proof}
Using Propositions~\ref {prop:twogenerators} and~\ref {prop:quadformexists},
we find a semiprimitive quadratic polynomial $A X^2 + B X + C$ modulo~$\mathfrak{n}$
representing the principally polarised ideal class.
It remains to apply a suitable matrix $\Mo$ as in
Proposition~\ref {prop:equivz} such that the resulting~$A'$ is totally
positive and coprime to~$\Ng$.
If $\pg$ is a prime ideal of~$\OKzero$ dividing~$\Ng$, we consider the
homogeneous form
$A' := A \delta^2 - B \gamma \delta + C \gamma^2$
in $\delta$ and $\gamma$ of \eqref {eq:equivform}.
Let
\begin {equation}
\label {eq:Aprimetop}
\begin{aligned}
M_\pg &= \id, & A' &= A & &\text {if }\pg\nmid A,\\
M_\pg &= \begin {pmatrix} 0 & -1 \\ 1 & \phantom{-}0 \end {pmatrix},
& A' &= C
& &\text{if } \pg \nmid C, \pg \mid A,\\
M_\pg &= \begin {pmatrix} 0 & -1 \\ 1 & \phantom{-}1
  \end {pmatrix},
&
A' &= A + C - B & & \text{otherwise (in which case $\pg \nmid B$)}.
\end{aligned}
\end {equation}
In all cases, we have $\pg \nmid A'$.

By Chinese remaindering, we obtain a matrix
$\Mo_{\Ng} \in \Sl_2 \left( \OKzero / \rad (\Ng) \right)$,
which can be lifted
to a matrix $\Mo \in \Sl_2 (\OKzero)$,
e.g.~by strong approximation~\cite[Appendix A.3]{garrett}.
Replace $z$ by $\Mo z$, so that $A$ gets replaced by $A'$,
which is coprime to~$\Ng$; and by Proposition~\ref {prop:equivz}
the total positivity of~$A$ carries over to~$A'$.
\end {proof}

\subsection {Class invariants}
\label {ssec:classinv}

We now have all ingredients at our disposal to state the first main result
of this article, which generalises the first statement in
\cite[Theorem~4, p.~331]{Schertz02} to CM fields of higher degree.
For a generalisation of the remainder of \cite[Theorem~4]{Schertz02},
see Theorem~\ref{th:conjugates} below.
Let
\[
\Gamma^0(N) = \left\{
\begin {pmatrix}
a & b \\ c & d
\end {pmatrix}
\in \Sp(\Z)
:
N\mid b\right\}.
\]
The first main theorem provides a sufficient criterion for a CM value of a function
invariant under $\Gamma^0 (N)$ to yield a class invariant,
cf.~Definition~\ref {def:classinv}.

\begin {theorem}
\label{th:maintheorem}
Suppose that $f \in \Fc_N$ is the quotient of two modular forms
with rational $q$-expansions and that it is invariant
under~$\Gamma^0 (N)$ for some positive $N\in\Z$.

Let $K_0$ be a totally real number field such that $\diff_{K_0}$ is
principal. Consider a polynomial $AX^2+BX+C\in\OKzero[X]$ of totally
negative discriminant $B^2 - 4 A C$ that is semiprimitive modulo~$N$
and such that $N \mid C$.
Let $z$ be a root of the polynomial and let $K = K_0 (z)$ be the
CM field generated by~$z$ over~$K_0$.
Let $\tau$ be obtained from $z$ as in
Propositions~\ref{prop:twogenerators}
and~\ref{prop:symplecticbasisgeneral}.

If $\tau$ is not a pole of~$f$, then $f (\tau)$ is a class invariant.
In other words, the value $f(\tau)$ lies in the Shimura class field of
the order and CM type corresponding to $z$ as in
Proposition~\ref{prop:twogenerators}.
\end{theorem}

The remainder of Section~\ref{ssec:classinv} is devoted to the proof
of Theorem~\ref{th:maintheorem}.
The main tool is Shimura's
reciprocity law, which describes the action of the Galois group of
$\HOPhiN / \Kr$ on CM values $f (\tau)$ by matrix actions
on the function~$f$.
For a commutative ring~$R$, let
\begin {equation}
\label {eq:GSp}
\GSp (R) = \left\{ \Mg \in \Mat_{2 g} (R) : \T {\Mg} \Jg \Mg = t \Jg
\text { for some } t \in R^\times \right\}
\end {equation}
with $\Jg$ as in~\eqref {eq:J}.
For $\Mg\in \GSp(R)$, write
\begin{equation}\label{eq:deft}
\T {\Mg} \Jg \Mg = t(\Mg) \Jg\qquad
\mbox{with}\qquad t(\Mg)\in R^\times.
\end{equation}
\begin{lemma}\label{lem:sympT}
	For $\Mg\in\GSp(R)$, both $\Mg^T$ and $\Mg^{-1}$ lie in $\GSp(R)$ as well.
	We have $t(\Mg^T)=t(\Mg)$ and $t(\Mg^{-1}) = t(\Mg)^{-1}$.
\end{lemma}
\begin{proof}
Moving $\T{\Mg}$ and $\Mg$ to the other side of the equality symbol in \eqref{eq:deft} yields
$\Jg = t(\Mg) \T{(\Mg^{-1})} \Jg \Mg^{-1}$, hence
$\Mg^{-1}\in \GSp(R)$ with $t(\Mg^{-1}) = t(\Mg)^{-1}$.
Subsequently inverting (and observing $\Jg^{-1} = -\Jg$) gives
$\Jg = t(\Mg)^{-1} \Mg \Jg \T{\Mg}$,
hence $\T{\Mg}\in \GSp(R)$ with $t(\T{\Mg}) = t(\Mg)$.
\end{proof}
It follows that $\GSp(R)$ is a subgroup of $\Gl_{2g}(R)$.
Let $\Sp (R)\subseteq \GSp(R)$ be the subgroup with $t = 1$ and
$\GSpplus (R)$ the subgroup with $t > 0$ for rings~$R$ where this definition
makes sense.
Notice that any matrix in $\GSp (R)$ can be written as
\[\begin {pmatrix} 1 & 0 \\ 0 & t \end {pmatrix} \Vg\qquad
\mbox{with} \qquad \Vg \in \Sp (R).\]

The action of $\Gamma = \Sp (\Z)$ on the Siegel space $\Hc_g$
as given in~\eqref {eq:siegelaction} extends to the
action $\tau \mapsto (a \tau + b)(c \tau + d)^{-1}$
of the full group $\GSpplus (\Q)$
and induces an action on Siegel modular functions by
$f^\Mg (\tau) = f (\Mg \tau)$ for $\Mg \in \GSpplus (\Q)$.
There is also a well-known action of $\GSp (\Z / N \Z)$ on~$\Fc_N$ as follows:
\begin {itemize}
\item
The action of a matrix in $\Sp (\Z / N \Z)$ is that of an
arbitrary lift to $\Sp (\Z)$.
\item
For $t \in (\Z / N \Z)^\times$, the matrix
$\begin {psmallmatrix} 1 & 0 \\ 0 & t \end {psmallmatrix}$ acts on the
$q$-coefficients of~$f$ as the Galois group element of
$\Q (\zeta_N) / \Q$ sending $\zeta_N$ to $\zeta_N^t$.
\end {itemize}

We use the following explicit formulation of Shimura reciprocity
as given in~\cite [Theorems 2.4 and~2.5]{Streng12}.

\begin {theorem}[Shimura's reciprocity law]
\label {th:reciprocity}
Let $\mathcal{O}$ be an order of $K$ and $\Phi$
a primitive CM type.
Let $(\bg,\xi)$ be a principally polarised
ideal for $(\Oc,\Phi)$ (cf.~Definition~\ref{def:polarisedclass}).
Let $\Sc$ be an $E_{\xi}$-symplectic basis of $\bg$
and let $\tau$ be the corresponding period matrix.
Then for any $f \in \Fc_N$ without a pole in~$\tau$, we have
$ f (\tau) \in \HOPhiN$.

Let $\sigma : \Cg_{\Oc, \Phi} (N) \to \Gal (\HOPhiN / \Kr)$
be the Artin map, cf.~\eqref{eq:artinmap}.
The Galois action on $f (\tau)$ is described as follows.

Let $F$ be the least positive integer such that $F\Oc_K\subset \Oc$.
For any $[\ag] \in \Cg_{\Oc, \Phi} (N)$ with $\ag \in \I (N F)$,
let $\Cc$ be a symplectic basis of $\cg = \NphirO (\ag)^{-1} \bg$
with respect to $E_{\Norm (\ag) \xi}$.
Let $\Mg\in\Gl_{2g}(\Q)$ be given by
$\Cc = \Sc \T {\Mg}$.
Then we have $\Mg \in \GSpplus (\Q)$ with $t=\Norm(\ag)^{-1}$,
and the reduction~$\mmod {\Mg}{N}$ exists and satisfies
$\mmod {\Mg}{N} \in \GSp (\Z / N \Z)$.
Then
\[\
f (\tau)^{\sigma ([\ag])}
= f^{(\mmod {\Mg}{N})^{-1}}  (\Mg \tau)
= f^{(\mmod {\Mg}{N})^{-1}}  (\tau'),
\]
where $\tau'$ is obtained from~$\Cc$
\end {theorem}

Shimura reciprocity links Galois actions to matrix actions on modular
functions and period matrices. In our setting the period matrices are
obtained from symplectic bases of ideals
as in~\eqref {eq:tauintermsofbasis}, written in terms of roots~$z$
of quadratic polynomials. We need to bring together the action of
symplectic matrices of dimension $2g \times 2g$ on period matrices
with the action of matrices of dimension $2 \times 2$ on elements
$z \in K$ as described in Proposition~\ref {prop:equivz}.
This can be done explicitly using the two bases
$\Bc_1 = (b_{1,1}, \ldots, b_{1,g})$ and
$\Bc_2 = (b_{2,1}, \ldots, b_{2, g})$ of $K_0$ over~$\Q$,
which we will interpret as matrices of dimension $1 \times g$ with
entries in $K_0$, occurring in our choice of the symplectic basis~$\Sc$
in Proposition~\ref{prop:symplecticbasisgeneral}.
The important property of these bases is that they are dual with respect
to the symmetric $\Q$-bilinear form
\[
K_0 \times K_0 \to \Q, \quad
(x, y) \mapsto L (x y)
\]
with $L$ a $\Q$-linear map, precisely,
\[
L : K_0 \to \Q, \quad
x \mapsto \Tr_{K_0 / \Q} \left( - \lambda^{-1} x \right);
\]
that is, $L (b_{1,i} b_{2,j}) = \delta_{i, j}$.
Given $\alpha\in K_0$ and $i,j\in\{1,2\}$,
denote by $[\alpha]_j^i\in\Mat_g(\Q)$
the transposed matrix of multiplication by $\alpha$ from $K_0$
with $\Q$-basis $\Bc_i$ to $K_0$ with $\Q$-basis $\Bc_j$,
that is,
\begin {equation}
\label {eq:blowup}
\alpha \T{\Bc_i} = \changewithcolumns{i}{j}{\alpha} \T{\Bc_j}.
\end {equation}
We obtain the following map from $\Mat_2 (K_0)$ to $\Mat_{2g} (\Q)$.

\begin {lemma}
\label {lem:M}
Let $\Mo = \begin {pmatrix} \alpha & \beta \\ \gamma & \delta \end {pmatrix}
\in \Mat_2 (K_0)$. Then
\begin {equation}
\label {eq:commuteblowup}
\T{\Mo} \begin {pmatrix} {\Bc_1} & 0 \\ 0 & {\Bc_2} \end {pmatrix}
=  \begin {pmatrix} {\Bc_1} & 0 \\ 0 & {\Bc_2} \end {pmatrix} \T{\otogleft{\Mo}}
\end {equation}
with
\begin {equation}
\label {eq:M}
\otogleft{\Mo} = \begin {pmatrix} \changewithcolumns{1}{1}{\alpha} & \changewithcolumns{1}{2}{\beta} \\
\changewithcolumns{2}{1}{\gamma} & \changewithcolumns{2}{2}{\delta} \end {pmatrix}
\in \Mat_{2g} (\Q).
\end {equation}
If $\Mo \in \Gl_2 (K_0)$ with $\det \Mo \in \Q$,
then $\otogleft{\Mo} \in \GSp (\Q)$ with $t(\otogleft{\Mo}) = \det \Mo$.
In particular, if $\Mo\in\Sl_2(\OKzero)$, then $\otogleft{\Mo}\in \Sp(\Z)$.
\end {lemma}

\begin {proof}
Formula~\eqref {eq:commuteblowup} follows by a direct computation from
the definition of~$[\cdot]_j^i$:
\[
\begin {pmatrix} \T{\Bc_1} & 0 \\ 0 & \T{\Bc_2} \end {pmatrix}\Mo
= \begin {pmatrix} \alpha \T{\Bc_1} & \beta \T{\Bc_1} \\ \gamma \T{\Bc_2} & \delta \T{\Bc_2} \end {pmatrix}
= \begin {pmatrix} \changewithcolumns{1}{1}{\alpha}\T{\Bc_1}  & \changewithcolumns{1}{2}{\beta} \T{\Bc_2 } \\
\changewithcolumns{2}{1}{\gamma}\T{\Bc_1}  & \changewithcolumns{2}{2}{\delta}\T{\Bc_2}  \end {pmatrix}
= \Mg \begin {pmatrix} \T{\Bc_1} & 0 \\ 0 & \T{\Bc_2} \end {pmatrix}
\]
with $\Mg = \otogleft{\Mo}$.
Concerning symplecticity of~$\Mg$, we use the map~$L$ and the two bases
$\Bc_1$ and $\Bc_2$ to define a bilinear form
$E : \Q^{2g} \times \Q^{2g} \to \Q$,
\begin{equation}
\label{eq:Efirst}
\begin {aligned}
E : (u,v)\ \longmapsto\ 
& L \left( \T{u}
           \T{\begin {pmatrix} \Bc_1 & 0 \\ 0 & \Bc_2 \end {pmatrix}}
           \begin {pmatrix} 0 & 1 \\ -1 & 0 \end {pmatrix}
           \begin {pmatrix} \Bc_1 & 0 \\ 0 & \Bc_2 \end {pmatrix} v
  \right) \\
& =
L \left( \T{u}
         \begin {pmatrix} 0 & \T {\Bc_1} \Bc_2 \\
                          - \T {\Bc_2} \Bc_1 & 0
         \end {pmatrix}
         v
  \right),
\end {aligned}
\end{equation}
in which the matrix between $\T{u}$ and $v$ is an element of
$\Mat_{2g} (K_0)$. Plugging in the unit vectors for $u$ and $v$ and using
that $\Bc_1$ and $\Bc_2$ are dual shows that~$E$ is a symplectic
form with the unit vectors as symplectic basis, so that in fact
$E (u, v) = \T{u} \Jg v$.
Plugging in $\T{\Mg} v$ for $v$ and $\T{\Mg} u$ for $u$ and
using~\eqref {eq:commuteblowup} changes the factor
\[
\begin {pmatrix} 0 & 1 \\ - 1 & 0 \end {pmatrix}
\text { into }
 {\Mo} \begin {pmatrix} 0 & 1 \\ - 1 & 0 \end {pmatrix} \T{\Mo}
= \det (\Mo) \, \begin {pmatrix} 0 & 1 \\ - 1 & 0 \end {pmatrix}.
\]
Since  $\det(\Mo)\in\Q$ can be taken outside the $\Q$-linear~$L$,
we get $E (\T{\Mg} u, \T{\Mg} v) = \det(\Mo) E(u,v)$.
In other words, we have $\Mg \Jg \T{\Mg} = \det(\Mo) \Jg$.
By Lemma~\ref{lem:sympT}, we get that $\Mg$ is symplectic
with $t(\Mg) = \det(\Mo)$.
\end {proof}

\begin{lemma}
\label{lem:MzMtau}
With $\tau(z)$ obtained from $z$ as in
Proposition~\ref{prop:symplecticbasisgeneral}, we have
\begin{equation}\label{eq:MzMtau}\tau(\Mo z) = \otogleft{\Mo} \tau(z),
\end{equation}
where $\otogleft{\Mo}$ is as in Lemma~\ref{lem:M}.
\end{lemma}

\begin{proof}
Take $z' = \Mo z$ and observe that Lemma~\ref{lem:M} gives
	\begin {align*}
	\Sc'
	:=&\ (z' \Bc_1 | \Bc_2)
	= (\gamma z + \delta)^{-1}
	\big( (\alpha z + \beta) \Bc_1 | (\gamma z + \delta) \Bc_2 \big) \\
	=&\ (\gamma z + \delta)^{-1}
	(z | 1) \, \T{\Mo} \, \begin {pmatrix} \Bc_1 & 0 \\ 0 & \Bc_2 \end {pmatrix}
	= (\gamma z + \delta)^{-1} \Sc \T{\otogleft{\Mo}}.
	\end {align*}
	The period matrix for the left hand side is $\tau(\Mo z)$.
	The period matrix for the right hand side is the same
	as the period matrix for $\Sc \T{\otogleft{\Mo}}$,
	which is $\otogleft{\Mo}\tau(z)$ by
	\cite[Lemma~4.7(b)]{Streng12}.
\end{proof}

\begin{remark}\label{rem:forgetfulmoduli}
Lemma~\ref{lem:MzMtau} tells us that $z\mapsto \tau(z)$ gives a well-defined
map from the set of $z$-values up to $\mathrm{SL}_2(\OKzero)$ to the set $\Sp(\Z)\backslash \Hc_g$.
In fact, this is exactly the natural map from the Hilbert moduli space to the Siegel moduli space restricted
to our values~$z$, which are CM points of the Hilbert moduli space.
\end{remark}

To treat polarised ideals only up to equivalence, we need to consider
multiplications by constants $\mu \in K^\times$; the following
lemma describes them explicitly as linear maps $K\rightarrow K$
with respect to different bases.
\newcommand{\oldMof}[1]{\T{\Mo_{#1}}}
\newcommand{\oldMmu}{\oldMof{\mu}}
\newcommand{\oldMMmu}{\T{\Mg_{\mu}}}
\newcommand{\newMof}[1]{\Mo_{#1}}
\newcommand{\newMmu}{\newMof{\mu}}
\newcommand{\newMMmu}{\Mg_{\mu}}
\newcommand{\newMMmuinv}{\Mg_{\mu^{-1}}}
\begin {lemma}
\label {lem:Mmu}
Let $A$, $B$, $C\in K_0$ and $z\in K$ be as in Theorem~\ref{th:maintheorem},
and let $\mu \in K^\times$.
Following~\eqref {eq:orderfromform} write
$\mu = \frac {\alpha A z + \beta}{d}$
with $\alpha \in \dg^{-1}$ for $\dg = \gcd (A, B, C)$,
$\beta \in \OKzero$ and $d \in \N$.

Then the matrices
\begin {equation}
\label {eq:Mmuprime}
\newMmu =
\frac{1}{d}
\begin {pmatrix}
\beta - \alpha B & -\alpha C \\
\alpha A & \beta
\end {pmatrix}
\in \Mat_2 \left( \frac {1}{d} \, \OKzero \right)
\end {equation}
and
\begin {equation}
\label {eq:Mmu}
\newMMmu = \otogleft{\newMmu} =
\frac{1}{d}
\begin {pmatrix}
\changewithcolumns{1}{1}{\beta-\alpha B} & \changewithcolumns{1}{2}{-\alpha C}\\
\changewithcolumns{2}{1}{\alpha A} & \changewithcolumns{2}{2}{\beta}
\end {pmatrix}
\in \Mat_{2g} \left( \frac {1}{d} \, \Z \right).
\end {equation}
satisfy
$\mu (z, 1) = (z, 1) \oldMmu$
and
$\mu\Sc = \Sc \oldMMmu$,
where $\Sc$ is the $\Q$-basis of~$K$ of
Proposition~\ref{prop:symplecticbasisgeneral}.
Furthermore $\newMMmu^{-1} = \newMMmuinv$,
where $\newMMmuinv$ is obtained from $\mu^{-1}$ in a
manner analogous to~\eqref {eq:Mmu}.
\end {lemma}

\begin {proof}
The identity $\mu (z, 1) = (z, 1) \oldMmu$ is obtained by direct computation
using $Az^2+Bz+C=0$.
Combining this with Lemma~\ref{lem:M}, we get
\begin {align*}
\mu \Sc &= \mu (z, 1)
\begin {pmatrix} \Bc_1 & 0 \\ 0 & \Bc_2 \end {pmatrix}
= (z, 1) \oldMmu
\begin {pmatrix} \Bc_1 & 0 \\ 0 & \Bc_2 \end {pmatrix} \\
&= (z, 1)
\begin {pmatrix} \Bc_1 & 0 \\ 0 & \Bc_2 \end {pmatrix}
\T\newMMmu
= \Sc \T\newMMmu.
\end {align*}
So $\newMMmu$ and $\newMMmuinv$ are the transposed matrices of
multiplication in $K$ by $\mu$ and $\mu^{-1}$ with respect to the same
$\Q$-basis, which implies that they are inverses of each other.
\end {proof}

\begin {proof}[Proof of Theorem~\ref{th:maintheorem}]
We use Theorem~\ref {th:reciprocity} to show that $f (\tau)$ is invariant
under
\[
\Gal \left( \HOPhiN / \HOPhione \right)
= \sigma \left(\frac{
\I (N F) \cap S_{\Oc, \Phi} (1) }{ S_{\Oc, \Phi} (N)}
\right).
\]
Let $\ag \in \I (N F) \cap S_{\Oc, \Phi} (1)$.
By the definition~\eqref {eq:ray}
of $S_{\Oc, \Phi} (1)$, there is some $\mu \in K^\times$
such that $\NphirO (\ag) = \mu \Oc$ and
$\Norm (\ag) = \mu \cconj \mu$.
As we took $\ag$ coprime to $N F$,
we have that
$\mu$ is coprime to~$N F$ with respect to~$\Oc$ by
condition~\point {\ref{point4}} of Definition~\ref {def:modstar}.

Let $\Sc$ be the symplectic basis of $\bg$ that gave rise to $\tau$.
Then $\Cc = \mu^{-1} \Sc$ is a symplectic basis
of $N_{\Phi^r,\Oc}(\ag)^{-1}\bg$.  We get $\Cc = \Sc \newMMmuinv^T$
in the notation of Lemma~\ref{lem:Mmu}.
Using Shimura's reciprocity law Theorem~\ref{th:reciprocity} we obtain
\[
f (\tau)^{\sigma ([\ag)]}
= f^{\left( {\mmod {(\newMMmuinv)}{N}} \right)^{-1}} (\tau')
= f^{{\mmod {(\newMMmu)}{N}}} (\tau')
= f^{{\mmod {(\newMMmu)}{N}}} (\tau),
\]
where we have used the last statement of Lemma~\ref {lem:Mmu} and that
$\tau'$, the period matrix obtained from $\Cc = \mu^{-1}\Sc$,
equals the period matrix obtained from~$\Sc$, which is $\tau$.
It remains only to show that $f^{{\mmod {(\newMMmu)}{N}}}=f$.

For this we write $\mu = \frac {\alpha A z + \beta}{d}$ as in
Lemma~\ref {lem:Mmu} with furthermore $d$ coprime to $NF$ (using
Definition~\ref {def:modstar}\point{\ref{point2}}
and Proposition~\ref{prop:quadformexists}).
Then by Theorem~\ref {th:reciprocity}, the matrix $\mmod {(\newMMmu)}{N}$,
as given by the four blocks of size $g \times g$
of~\eqref {eq:Mmu}, is an element of $\GSp (\Z / N \Z)$ with
$t = \Norm (\ag) = \mu \overline \mu \pmod N$
coprime to~$N$.
So ${\mmod {(\newMMmu)}{N}}$ is the product of the matrix
$
\begin {psmallmatrix}
1 & 0 \\ 0 & t
\end {psmallmatrix}
\in \Mat_{2g} (\Z / N Z),
$
under which the function~$f$ is invariant as a quotient of forms with
rational $q$-expansions,
and a matrix in $\Sp (\Z / N \Z)$ for which the upper right block equals
$\changewithcolumns{1}{2}{- \alpha C} \bmod N$.
Since $\alpha \in \dg^{-1}$ we have $\alpha C \in \OKzero$.
With $\dg = \gcd(A,B,C)$ coprime to~$N$, the element $\alpha\in \dg^{-1}$
has non-negative
valuation in all primes of~$K_0$ dividing $N$, so that $N \mid C$ implies
$\alpha C \in N\OKzero$. Thus the upper right block vanishes modulo~$N$,
and the matrix in $\Sp (\Z / N \Z)$ can be lifted to a matrix in
$\Gamma^0 (N)$, under which the function~$f$ is also invariant by
assumption.
\end {proof}

\begin{remark}
	\label{rem:compareLauterYang}
	Our map $\Mo\mapsto \otogleft{\Mo}$ generalises the map~$\phi$ of \cite[(3.4)]{LauterYang}.
	They have $K_0 = \Q(\sqrt{D})\subset \R$ with
	$D\equiv 1\pmod{4}$ prime and $\sqrt{D}>0$,
	let $\sigma$ be the non-trivial automorphism of $K_0$,
	write  $\OKzero = \Z e_1 + \Z e_2$ and take $\epsilon \in \OKzero^\times$ with $\epsilon\sigma(\epsilon) = -1$ and $\epsilon>0$.
	Observe that $e_1\sigma(e_2) - e_2\sigma(e_1) = s\sqrt{D}$
	for some $s\in\{\pm 1\}$.
	Then loc.~cit.\ defines $\phi(\Mo) = S\Mo^* S^{-1}$, where
	$$	\Mo^* = \begin{pmatrix} \alpha & 0 & \beta & 0\\
		0 & \sigma(\alpha) & 0 & \sigma(\beta)\\
		\gamma & 0 & \delta & 0\\
		0 & \sigma(\gamma) & 0 & \sigma(\delta)\end{pmatrix}\quad \mbox{and}\quad
	S = \begin{pmatrix}
		e_1 & \sigma(e_1) & 0 & 0\\
		e_2 & \sigma(e_2) & 0 & 0\\
		0 & 0 & \frac{s\sigma(e_2)}{\epsilon} & \frac{s e_2}{\sigma(\epsilon)}\\
		0 & 0 & \frac{-s\sigma(e_1)}{\epsilon} & \frac{-s e_1}{\sigma(\epsilon)}
	\end{pmatrix}.
	$$
	Now choose $\Bc_1 = (e_1, e_2)$, which has
	trace dual $s\sqrt{D}^{-1}(\sigma(e_2), -\sigma(e_1))$.
	Choosing $\lambda = -\sqrt{D}/\epsilon$ in Proposition~\ref{prop:symplecticbasisgeneral},
	we obtain $\Bc_2$ as $-\lambda$ times the trace dual of $\Bc_1$, which is
	$\Bc_2 = (\frac{s}{\epsilon}\sigma(e_2), -\frac{s}{\epsilon}\sigma(e_1))$.
	In particular, the first and third column of the matrix identity $S\Mo^* = \phi(\Mo) S $
	read
	$$	S \begin{pmatrix} \alpha & \beta\\
		0 & 0\\
		\gamma & \delta\\
		0 & 0\end{pmatrix} = \phi(\Mo) \begin{pmatrix}
		e_1 & 0 \\
		e_2 & 0\\
		0 & \frac{s\sigma(e_2)}{\epsilon} \\
		0 & \frac{-s\sigma(e_1)}{\epsilon}
	\end{pmatrix},
	$$
	which is exactly 	
	$$ \begin{pmatrix} \T{\Bc_1} & 0 \\ 0 & \T{\Bc_2} \end{pmatrix} \Mo
	= \phi(\Mo) \begin{pmatrix} \T{\Bc_1} & 0 \\ 0 & \T{\Bc_2} \end{pmatrix}.$$
	As \cite[(3.4)]{LauterYang} gives that $\phi(\Mo)$ has rational coefficients,
	this implies that it equals~$\otogleft{\Mo}$,
	so that indeed the map $\Mo\mapsto \otogleft{\Mo}$ generalises~$\phi$.
\end{remark}

\subsection {Existence of quadratic polynomials with
\texorpdfstring {$N \mid C$}{N dividing C}}
\label {ssec:existenceform}

We would like to apply Theorem~\ref{th:maintheorem}
to arbitrary orders $\Oc$ and integers~$N$.
The requirements of the theorem are twofold:
On the one hand, the function needs to be invariant under some $\Gamma^0 (N)$.
Such functions are plentiful, and we provide some interesting families
of examples in~\S\ref {sec:functions}.
On the other hand, we need the existence of a suitable quadratic polynomial;
using the terminology of Definitions~\ref {def:polarisedclass}
and~\ref {def:represents}, we need the existence of a polarised
ideal class $T$ for~$(\Oc, \Phi)$ that is represented by a
quadratic polynomial $Q = A X^2 + B X + C$ satisfying
\begin {equation}
\label{eq:requirementsmainthm}
Q \text { is semiprimitive modulo } \Ng \text { and } \Ng \mid C
\end {equation}
for $\Ng = N \OKzero$.
The following theorem gives a necessary and sufficient criterion for
the existence of such a polynomial in the case that $\Ng$ is prime to the
conductor, which includes the particularly important case $\Oc = \OK$.
Since the proof is identical, we formulate it directly for the case of a
general ideal~$\Ng$ of~$\OKzero$, although later applications will only
need the case $\Ng = N \OKzero$.
The result assumes the technical condition that a polarised ideal class for
$(\Oc, \Phi)$ exists, but otherwise the question of computing a class
polynomial would be moot.

\begin {theorem}
\label {th:n}
Let $\Ng$ be a non-zero ideal of $\OKzero$,
assume that $\diff_{K_0} = \diffgen \OKzero$ and that
$\Oc\subseteq K$ is an order
of conductor~$\fg$ coprime to $\Ng$ and containing~$\OKzero$,
and let $F \Z = \fg \cap \Z$.
Suppose that there exists a principally polarised ideal class for
$(\Oc,\Phi)$.
Then the following are equivalent:
\begin{enumerate}
\item Every prime ideal of $\OKzero$ dividing $\Ng$
is either split in $\OK$, or it is ramified
and divides~$\Ng$ with multiplicity~$1$.
\item Every principally polarised ideal class for $(\Oc,\Phi)$
is represented (as in Definition~\ref{def:represents})
by a polynomial satisfying~\eqref{eq:requirementsmainthm}
with
\[
\gcd (A, F \Ng) = 1 \text { and } \gcd(\Ng, \Ng^{-1} C) = 1.
\]
\item There exists a principally polarised ideal
class for $(\Oc,\Phi)$ that is represented by a polynomial
satisfying~\eqref{eq:requirementsmainthm}.
\end{enumerate}
Furthermore, if \point {3} holds and $\OKzero$ has
narrow class number~$1$,
then
\begin{itemize}
\item[\point {3'}] the assertion of \point {3} holds with $A = 1$;
\item[\point {3''}] the assertion of \point {3} holds with $C = \nu$,
   where $\Ng = \nu \OKzero$.
\end{itemize}
\end {theorem}

We will use the following special case of the Kummer-Dedekind
theorem in the proof.
\begin {lemma}
\label {lem:dedekind}
Let $\Oc \subseteq K$ be an order of conductor~$\fg$
and containing $\OKzero$,
assume $\diff_{K_0} = \diffgen \OKzero$, and
let $\pg$ be a prime ideal of~$\OKzero$ not dividing~$\fg$, and let $z \in K$
be a root of a quadratic polynomial $A X^2 + B X + C$ as in
Proposition~\ref {prop:quadformexists} with
$\pg \nmid \dg = \gcd (A, B, C)$ and such that $z \OKzero + \OKzero$ has
multiplier ring~$\Oc$.
Write $U (X) = X^2 + B X + A C \in \OKzero [X]$ with root $\theta = A z$,
and let $\widetilde U$ be the reduction of $U$ modulo~$\pg$.

Then the splitting behaviour of~$\pg$ in~$\Oc$ is governed by the
factorisation of $\widetilde U$ in $(\OKzero / \pg) [X]$ as follows.
If $\widetilde U = \prod_i \widetilde U_i^{e_i}$ with
monic $\widetilde {U_i}$ and $U_i$ is an arbitrary monic lift
of $\widetilde U_i$ to $\OKzero [X]$, then the ideals above $\pg$
are given by the $\Pg_i = \pg \Oc + U_i (\theta) \Oc$ of residue field
degree~$f_i = \deg \widetilde U_i$ and ramification index~$e_i$.

Moreover, if $e_1=2$, then the remainder of $U$ upon division by $U_1$,
which is an element of~$\OKzero$, is not divisible by~$\pg^2$.
In particular, $U$ has no root modulo~$\pg^2$.
\end {lemma}

\begin {proof}
Notice that since $\pg$ is coprime to~$\fg$, its splitting in~$\Oc$ is the
same as in~$\OK$.
By \eqref {eq:orderfromform},
we have $\theta = A z \in \Oc$ and
\[
\Oc = \dg^{-1} \theta + \OKzero,
\]
which implies that the conductor of $\OKzero [\theta]$ divides $\fg \dg$.
As $\pg \nmid \fg \dg$, the Kummer-Dedekind criterion for relative number field
extensions gives the statements of the lemma;
see \cite[Proposition~2.3.9]{Cohen00}
and~\cite[Theorem 8.2]{stevenhagen-number-rings}.

As the first reference does not include the final statement
and the second reference states the results only for $\OKzero=\Z$,
we carry out the proof of the final statement.
Write $U (X) = q (X) U_1 (X) + r$ with $q \in \OKzero [X]$ monic and
linear and $r \in \OKzero$. From $v_{\Pg_1} (U_1 (\theta)) = 1$
and $U (\theta) = 0$ we deduce $\Pg_1 \mid r$, which is equivalent to
$\pg \mid r$ since $r \in \OKzero$. This implies
$\widetilde q = \widetilde {U_1}$, so that $v_{\Pg_1} (q (\theta)) = 1$
and
$v_\pg (r) = \frac {1}{2} v_{\Pg_1} (r)
= \frac {1}{2} v_{\Pg_1} (q (\theta) U_1 (\theta)) = 1$.
If $U$ had a root modulo~$\pg^2$, then we could choose without loss of
generality $U_1$ such that it would have this root modulo~$\pg^2$,
which would imply the contradiction $\pg^2 \mid r$.
\end {proof}

\begin {proof}[Proof of Theorem~\ref{th:n}]
The implication \point{2}$\Rightarrow$\point{3} is trivial under the
assumption that some polarised ideal class exists for~$(\Oc, \Phi)$.

We start with  \point{3}$\Rightarrow$\point{1}.
Assume that $z$ is the root of a polynomial
satisfying~\eqref {eq:requirementsmainthm}
and that $(\bg, \xi)$ is the associated principally polarised ideal
as in Proposition~\ref {prop:twogenerators}.
Every prime $\pg\mid \Ng$ satisfies $\pg \mid C$.
In particular, using the notation of Lemma~\ref {lem:dedekind},
$\widetilde{U}$ is reducible, so the prime $\pg$
is not inert.
If $\pg$ is ramified, that is $\pg \mid B$,
then we have $\Pg_1 = \pg \OK + \theta \OK$
with $v_{\Pg_1} (\pg \OK) = 2$, hence
$v_{\Pg_1} (\theta) = 1$ and
$v_\pg (\Ng) \leq v_\pg (C) = v_\pg (AC)
= v_\pg (\Norm_{K / K_0} (\theta)) = 1$.

Now we prove \point{1}$\Rightarrow$\point{2}.
Let $T$ be a principally polarised ideal class for $(\Oc, \Phi)$; by
Proposition~\ref {prop:Acoprime}, it can be represented by a
quadratic polynomial $A X^2 + BX+C$ with $A \gg 0$ and
$\gcd (A, F \Ng) = 1$.
We show how to modify~$z$ such that furthermore
$\Ng \mid C$.
Let $\pg$ be a prime dividing~$\Ng$. As it is coprime to~$\fg$ and~$A$
and split or ramified, the polynomial $U(X) = X^2 + BX + AC$
of Lemma~\ref {lem:dedekind} has a root in $\OKzero / \pg$.
If $\pg$ is split, then this root is simple, so we may Hensel lift it to a
root modulo an arbitrary power of~$\pg$. The Chinese remainder theorem
allows us to combine the roots into an element $r \in \OKzero$ that is a root modulo~$\Ng$.
As $A$ is coprime to~$\Ng$, we may furthermore assume that $A \mid r$.
Let $\theta' = \theta - r$; its minimal polynomial is
$U' = U (X + r) = X^2 + B' X + AC'$ with $B' = B + 2r$ and $C' = U(r)/A \in \Ng$.
So $z' = \theta' / A$ is a root of the polynomial
$A X^2 + B' X + C' \in \OKzero$
satisfying~\eqref {eq:requirementsmainthm}.
The polynomial is obtained by the $\Sl_2(\OKzero)$-transformation
$z' = z-r/A$, where $A\mid r$, so by Proposition~\ref {prop:equivz}
it represents the same principally polarised ideal class.

We now refine this argument so as to obtain
$\gcd(\Ng,\Ng^{-1}C)=1$, that is,
all primes $\pg$ dividing~$\Ng$
satisfy $v_{\pg}(\Ng)=v_{\pg}(C)$.
Given a prime $\pg \mid \Ng$, let $e=v_{\pg}(\Ng)$.
If $\pg$ splits, then there are \emph{unique}
Hensel lifts of each of the two distinct roots
of~$\widetilde U$ modulo~$\pg$ to a root modulo~$\pg^e$ and a root
modulo~$\pg^{e+1}$.
As each element of $\OKzero / \pg^e$ has $\Norm_{K_0/\Q} (\pg) \geq 2$
different lifts to an element of $\OKzero / \pg^{e+1}$, we may choose a root
modulo~$\pg^e$ that is not a root modulo~$\pg^{e+1}$. In this way, the final
$C'$ is divisible by~$\pg^e$, but not by $\pg^{e+1}$.
If $\pg$ ramifies in $K/K_0$, then $e=1$
and by Lemma~\ref{lem:dedekind} the quadratic polynomial
has no root modulo~$\pg^2$, so that
$v_{\pg}(\Ng)=v_{\pg}(C)$ is automatically true.

Assume now that $K_0$ has narrow
class number~$1$ and that \point{3} holds.
It remains to prove that \point{3'} and \point{3''} hold.
We start with the proof of~\point{3''}.
Write $\Ng=\nu\OKzero$ with~$C / \nu$ totally positive.
We have already
reached $\nu \mid C$ and $\gcd (A, \nu) = \gcd (C / \nu, \nu) =~1$.
Also without loss of generality we assume $\dg=\gcd(A,B,C)=~1$,
even with $A \gg 0$, by the requirement on the narrow class number
and $\gcd (\dg, \Ng) = 1$.

Then let $z' = z \nu / C$,
$A' = A \frac {C}{\nu} \gg 0$,
$C' = \nu$ and $B' = B$.
We still have $\gcd (A', \nu) = 1$
and hence $\dg'=\gcd(A',B',C')=1=\dg$.
As we also have $A'z'=Az$, we find that $z'\OKzero+\OKzero$
has the same endomorphism ring $\Oc$ as $z\OKzero+\OKzero$
by \eqref {eq:orderfromform},
and since $\nu/C \gg 0$ we find that it has the same CM type.
This finishes the proof of~\point{3''}.

The proof of~\point{3'} is exactly the same, but
with $z'=Az$, $A'=1$ and $C'=AC$.
\end {proof}

\section{\texorpdfstring{$N$}{N}-systems}
\label {sec:nsystems}

The main Theorem~\ref {th:maintheorem} of the preceding section provides a
convenient and very generic way of obtaining class invariants in the sense
of Definition~\ref {def:classinv}. For computing them algebraically, we
need a handle on their characteristic polynomials (see also
Definition~\ref {def:classinv}); otherwise said, we need to explicitly
describe their Galois conjugates.

A tool for doing so, introduced for $g = 1$ in \cite {Schertz02},
are \emph{$N$-systems}, consisting of quadratic polynomials representing (in the sense of
Definition~\ref {def:represents}) the equivalence classes of principally
polarised ideals and satisfying certain congruence conditions modulo~$N$.
We generalise this notion to arbitrary~$g$ in \S\ref {ssec:nsystemstheory}
and prove that an $N$-system describes a complete set of Galois conjugates
of class invariants. Then in \S\ref {ssec:nsystemspractice} we show that
$N$-systems always exist and provide an algorithm to compute them explicitly.

Throughout this section we assume as before that the order $\Oc$ of
conductor~$\fg$ is closed under complex conjugation and contains~$\OKzero$,
that $\diff_{K_0} = \diffgen\OKzero$ and that $F \Z = \fg \cap \Z$.

\subsection {Galois conjugates from \texorpdfstring {$N$}{N}-systems}
\label {ssec:nsystemstheory}

As seen in \S\ref {ssec:ray}, the CM class group of level~$1$,
$\Cg_{\Oc, \Phi} (1)$, acts freely on the set $\Tc_{\Oc, \Phi}$ of
principally polarised ideal classes via~\eqref {eq:cmgroupaction}.
Let $\Tc = \{T_1,\ldots,T_h\} \subseteq \Tc_{\Oc, \Phi}$ with
$h = |\Cg_{\Oc, \Phi} (1)|$ be one orbit under this action.
In the light of Proposition~\ref {prop:Acoprime}, the classes in~$\Tc$
may be
represented by quadratic polynomials $A_i X^2 + B_i X + C_i \in \OKzero [X]$,
the roots~$z_i$ of which determine period matrices~$\tau_i$ as in
Proposition~\ref {prop:symplecticbasisgeneral} and
Corollary~\ref {cor:symplecticbasisg2}.

If $f$ is a Siegel modular function of level~$N = 1$, then Shimura
reciprocity implies that for any choice of quadratic polynomials, the
$f (\tau_i)$ form a Galois orbit and are thus the roots of a class
polynomial as in~\eqref {eq:classpol}.
For general~$N$, each $T_i$ can be given by
$[S_{\Oc, \Phi} (1) : S_{\Oc, \Phi} (N)]$ representatives that are
inequivalent under the action of $S_{\Oc, \Phi} (N)$, and which in general
yield different values of~$f$;
so some work is needed to write down a consistent set of quadratic
polynomials leading to a Galois invariant set of CM values.
Before giving a solution by imposing congruence conditions modulo~$N$,
we address another problem arising for $g > 1$. When $g = 1$, it is
natural to choose \emph{primitive} representatives, satisfying $A_i > 0$
and $\gcd (A_i, B_i, C_i)=1$. In the general case, this rigidity is not
possible unless the narrow class number of~$\OKzero$ is~$1$.
So far we worked around this by using the notion of semiprimitivity modulo
the ideal $\Ng = N \OKzero$, see Definition~\ref {def:semiprimitive} and
Proposition~\ref {prop:quadformexists}, which imposes that the greatest
common divisor of the coefficients is coprime to~$\Ng$. This is compatible
with multiplying $A_i$, $B_i$ and $C_i$ by the same element of~$\OKzero$
coprime to~$\Ng$ without changing~$z_i$ and~$\tau_i$, which turns out to be
undesirable. So we impose an additional notion of compatibility between the
quadratic polynomials.

Again we state results in terms of an arbitrary non-zero ideal
$\Ng$ of $\OKzero$ when their proofs are rigorously identical;
this will be useful in future work using Hilbert modular forms.
However, for the applications in the remainder of this article only the case
that $\Ng$ is generated by a rational integer will be needed.

\begin{definition}
\label {def:equiprimitive}
Let $\Ng$ be an integral ideal of $\OKzero$.
A pair of quadratic polynomials $A_1 X^2+B_1 X+C_1$ and
$A_2 X^2+B_2 X+C_2\in\OKzero [X]$
is \emph{equiprimitive modulo $\Ng$} if
both are semiprimitive modulo $\Ng$ and
their discriminants $D_1 = B_1^2 - 4 A_1 C_1$ and $D_2 = B_2^2 - 4 A_2 C_2$
are equal.
\end{definition}

The following lemma shows that equiprimitivity forces the greatest common
divisors of the coefficients to all be the same over the set of quadratic
polynomials.
\begin{lemma}
\label{lem:sqrtd}
Assume $\Oc \supseteq \OKzero$ and $\diff_{K_0} = \diffgen \OKzero$.
Let $A_i X^2 + B_i X + C_i \in \OKzero [X]$ with roots $z_i$
for $i \in \{ 1, 2 \}$
represent two classes of principally polarised ideals for $(\Oc, \Phi)$.
Write
$\dg_i = \gcd (A_i, B_i, C_i)$,
$\delta_i = 2 A_i z_i + B_i$
and $\epsilon = \delta_1 \delta_2^{-1}$.
Then
\begin{equation}
\label{eq:newepsilon}
\epsilon \OKzero = \dg_1 \dg_2^{-1}.
\end{equation}
If the two quadratic polynomials are semiprimitive modulo~$\Ng$,
then $\epsilon$ is coprime to~$\Ng$ and totally positive.

If the two quadratic polynomials are equiprimitive modulo~$\Ng$,
then $\epsilon = 1$, that is, $\delta_1=\delta_2$ and $\dg_1 = \dg_2$.
\end {lemma}

\begin {proof}
Notice that $\delta_i^2 = D_i = B_i^2 - 4 A_i C_i$, so that $\delta_i$
is a square root of the discriminant~$D_i$.
Notice also that $\epsilon = \frac {A_1 \xi_2}{A_2 \xi_1}$ for
$\xi_i = \left( (z_i - \cconj z_i) \diffgen \right)^{-1}$, which are purely
imaginary; so $\epsilon$ is an
element of~$K_0$.

From~\eqref {eq:orderfromform} we have the two expressions for
$\Oc$ as $\Oc = \dg_i^{-1} A_i z_i + \OKzero$, leading to
\[
2 \Oc + \OKzero
= \dg_i^{-1} 2 A_i z_i + \OKzero
= \dg_i^{-1} (2 A_i z_i + B_i) + \OKzero
= \dg_i^{-1} \delta_i + \OKzero,
\]
so that
\[
\dg_1^{-1} \delta_1 + \OKzero = \dg_2^{-1} \delta_2 + \OKzero.
\]
Since $\OKzero$ and the $\dg_i$ are real and the $\delta_i$ are purely
imaginary, we may ``compare imaginary parts'' and find the desired
equality~\eqref{eq:newepsilon}.

In the semiprimitive case, by definition the $\dg_i$ are coprime to~$\Ng$ and
the $A_i$ are totally positive. So $\epsilon \OKzero = \dg_1 \dg_2^{-1}$ is
also coprime to~$\Ng$. Moreover,
the signs of the two real embeddings of $\epsilon$ are those of the embeddings
of $\xi_2 \xi_1^{-1}$ under the CM type, and since the $\xi_i$ have positive
purely imaginary embeddings, their quotient is totally positive.

In the equiprimitive case, the element $\epsilon$ is
a totally positive square root of $D_1 / D_2 = 1$, so $\epsilon = 1$,
which means $\delta_1 = \delta_2$ and $\dg_1 = \dg_2$.
\end {proof}

We postpone further discussion of the properties of equiprimitive polynomials
to \S\ref {ssec:nsystemspractice} in favour of stating the definition of
$N$-systems and the main result on Galois conjugates.

\begin{definition}
\label{def:nsystem}
Let $\Ng$ be a non-zero ideal of $\OKzero$,
and let
$\Qc = \{ Q_1, \ldots, Q_h \}$ with
$Q_i = A_i X^2 + B_i X + C_i \in \OKzero [X]$
be a set of polynomials representing an orbit of principally polarised
ideal classes under the action of $\Cg_{\Oc, \Phi} (1)$.
We call $\Qc$ an \emph{$\Ng$-system for $(\Oc,\Phi)$}
if it consists of equiprimitive polynomials modulo $F\Ng$ that satisfy
\begin{enumerate}
\item\label{item:nsystem2} $\gcd (A_i, F \Ng)=1$,
\item\label{item:nsystem4} $B_i\equiv B_j \pmod {2 \Ng}$
for all~$i$ and~$j$.
\end{enumerate}
If $\Ng = N \OKzero$ for some $N\in\Z_{>0}$, then we call
$\Qc$ an \emph {$N$-system}.
\end{definition}

In the case $g=1$, the action of the Galois group on the ideal class group
is transitive, and every $N$-system as in \cite[p.~329]{Schertz02}
is also an $N$-system in our sense.
Compared to \cite {Schertz02} we have added the condition that
$\gcd (A, F) = 1$, which simplifies some proofs without having practical
implications, as it is always possible by Theorem~\ref {th:nsystems} and
Algorithm~\ref {alg:nsystem} to satisfy the stricter condition.
The following is a generalisation
of Schertz~\cite[Theorem~4, pp.~331--332]{Schertz02}
from the case $g=1$.

\begin{theorem}
\label{th:conjugates}
Under the hypotheses of Theorem~\ref{th:maintheorem},
let $\tau_1, \ldots, \tau_h$ be the period matrices obtained
as in Proposition~\ref {prop:symplecticbasisgeneral} and
Corollary~\ref {cor:symplecticbasisg2} from
an $N$-system given as $\{ Q_1,\ldots, Q_h \}$
such that $N \mid C_1$ and $\tau_1$ is not a pole of~$f$.

Then $f(\tau_1)$ is a class invariant,
and the set of Galois conjugates of $f (\tau_1)$ over~$K^r$
is exactly $\{ f (\tau_1), \ldots, f (\tau_h) \}$.

More precisely, let~$(\bg_i, \xi_i)$ be the polarised ideal associated
to~$Q_i$. Then there is $\ag_i \in \I (N F)$ and $\mu_i\in K^\times$
such that
$\bg_i = \mu_i^{-1} \NphirO (\ag_i)^{-1} \bg_1$ and
$\xi_i = \mu_i\cconj{\mu_i} \Norm (\ag_i) \xi_1$.
For any such $\ag_i$, we have
\begin {equation}
\label {eq:artin}
f (\tau_1)^{\sigma ([\ag_i])} = f (\tau_i).
\end {equation}
\end{theorem}

\begin{proof}
We know from Theorem~\ref{th:maintheorem} that $f(\tau_1)$ is a class
invariant.
Equiprimitivity and the other properties of an $N$-system imply that
all the $C_i$ are divisible by $N$; so actually all the $f (\tau_i)$ are
class invariants. It remains to show that they are Galois conjugates
as given by~\eqref {eq:artin}.

The $\ag_i$ and $\mu_i$ exist because by definition the polarised ideal
classes derived from an $N$-system form an orbit under the
action~\eqref {eq:cmgroupaction} of the CM class group $\Cg_{\Oc, \Phi} (1)$;
here a priori $\ag_i \in \I (F)$, but (as usual in class field theory)
we may take a representative satisfying the additional coprimality condition
$\ag_i \in \I (N F)$.\label{page:proofofthconjugates}

The action of $\sigma ([\ag_i])$ is computed using
Theorem~\ref {th:reciprocity} as follows.
With the notations of Theorem~\ref{th:reciprocity} and
Proposition~\ref {prop:symplecticbasisgeneral}, we have
$\cg = \NphirO (\ag_i)^{-1} \bg_1 = \mu_i \bg_i$,
\[
\Sc = (z_1, 1)
\begin {pmatrix} \Bc_1 & 0 \\ 0 & \Bc_2 \end {pmatrix}
\text { and }
\Cc = \mu_i (z_i, 1)
\begin {pmatrix} \Bc_1 & 0 \\ 0 & \Bc_2 \end {pmatrix}.
\]
By equiprimitivity of an $N$-system and Lemma~\ref {lem:sqrtd}, we have
\[
2 A_i z_i + B_i = 2 A_1 z_1 + B_1,
\text { so }
(z_i, 1) = (z_1, 1) \T{\Mo}
\text { with }
\Mo = \begin {pmatrix}
\frac {A_1}{A_i} & \frac {B_1 - B_i}{2 A_i} \\
0 & 1
\end {pmatrix}.
\]
Let $\newMof{\mu_i}$ be as in Lemma~\ref{lem:Mmu} for $z=z_1$, that is,
such that $\mu_i (z_1,1) = (z_1,1) \oldMof{\mu_i}$.
Then
\begin {align*}
\Cc
&= \mu_i (z_1, 1) \T{\Mo}
\begin {pmatrix} \Bc_1 & 0 \\ 0 & \Bc_2 \end {pmatrix}
= (z_1, 1) \oldMof{\mu_i} \T{\Mo}
\begin {pmatrix} \Bc_1 & 0 \\ 0 & \Bc_2 \end {pmatrix}\\
&= (z_1, 1) \T{(\Mo \newMof{\mu_i})}
\begin {pmatrix} \Bc_1 & 0 \\ 0 & \Bc_2 \end {pmatrix}\\
&= \Sc  \T{\Mg} \quad\mbox{with}\quad \Mg = \otogleft{\Mo\newMof{\mu_i}}\quad\mbox{as in Lemma~\ref{lem:M}}.
\end {align*}
By Shimura's reciprocity law (Theorem~\ref{th:reciprocity}) we get
\begin{equation}\label{eq:resultshimrecipinproof}
f(\tau_1)^{\sigma([\ag_i])} =  f^{\mmod {\Mg}{N}^{-1}}(\tau'),
\end{equation}
where $\tau'$ is obtained from the basis $\Sc \T{\Mg}$, hence $\tau' = \tau_i$.
As in the proof of Theorem~\ref{th:maintheorem},
we will show $f^{\mmod {\Mg}{N}^{-1}}=f$ by looking at the upper right
entry of~$\Mg$.

Write $\theta_i = A_i z_i \in \Oc$ and $\dg_i = \gcd (A_i, B_i, C_i)$
and notice that $A_i \bg_i = \theta_i \OKzero + A_i \OKzero$
is an integral ideal of~$\Oc$
and that $A_i \cconj \bg_i$ is an integral ideal of~$\Oc = \cconj \Oc$.
Then
\begin {equation}
\label {eq:normb}
\begin {split}
(A_i \bg_i) (A_i \cconj \bg_i)
& = A_i (A_i z_i \OKzero + A_i \OKzero)(\cconj z_i \OKzero + \OKzero) \\
& = A_i (A_i z_i \cconj z_i \OKzero + A_i (z_i + \cconj z_i) \OKzero
    + A_i \OKzero + \theta_i \OKzero) \\
& = A_i (\dg_i + \theta_i \OKzero)
= A_i \dg_i \Oc
\text { by \eqref {eq:orderfromform}}.
\end {split}
\end {equation}

As $A_i$ is coprime to $N F \OKzero$, this shows that all the
$\bg_i$ (including $\bg_1$) are coprime to $N F \Oc$; with
$N_{\Phir, \Oc}(\ag_i)$ being coprime to $N F \Oc$, this implies
by Definition~\ref {def:modstar}\point {\ref{point4}}
that $\mu_i$ is coprime to $N F$ with respect to $\Oc$.
We may write it as in Lemma~\ref {lem:Mmu} as
$\mu_i = \frac {\alpha_i A_1 z_1 + \beta_i}{d_i}$ with
$\alpha_i \in \dg_1^{-1}$, $\beta_i \in \OKzero$ and,
by Definition~\ref {def:modstar}\point {\ref{point2}}, with a denominator
$d_i \in \Z$ that is coprime to~$N F$.
As $C_1\in N\OKzero$ by assumption,
we then see from~\eqref {eq:Mmuprime}
that the top right entry of $\newMof{\mu_i}$
is divisible by~$N$
(in the sense that its valuation in every prime ideal
$\pg \mid N \OKzero$ is at least $v_{\pg} (N \OKzero)$).

By properties~\point{\ref{item:nsystem2}}
and~\point {\ref{item:nsystem4}} of an $N$-system
we find that the top right entry of $\Mo$ is also
divisible by~$N$,
and hence the same holds for the product of the $NF$-integral matrices
$\Mo$ and $\newMof{\mu_i}$.

The matrix ${\Mg} = \otogleft{\Mo \newMof{\mu_i}}$ is obtained from
$\Mo \newMof{\mu_i}$ as in~\eqref {eq:M}
by replacing elements of $K_0$ by their $g \times g$-matrices
with respect to $\Z$-bases of $\OKzero$.
In particular, if an element of $K_0$ is $N$-integral, then
so are the entries of the corresponding $g\times g$ block.
And if an element of $K_0$ is divisible by~$N$, then
so are the entries of the corresponding $g\times g$ block.
So the transposed matrix $\Mg$ is $N$-integral with upper right block
divisible by~$N$.
We conclude
$f^{\mmod {\Mg}{N}^{-1}} = f$,
which finishes the proof.
\end{proof}

\subsection {Existence and computation of \texorpdfstring {$N$}{N}-systems}
\label {ssec:nsystemspractice}

We show that an $\Ng$-system in the sense of Definition~\ref {def:nsystem}
always exists by describing an algorithm to transform any set of
polynomials representing an orbit of principally polarised ideal classes
into an $\Ng$-system.
The following is a generalisation of
Schertz~\cite[Proposition~3, pp.~335--336]{Schertz02}.

\begin{theorem}
\label{th:nsystems}
Let $\Oc$ be an order in a CM field~$K$ that is closed under complex
conjugation and contains~$\OKzero$, assume that
$\diff_{K_0} = \diffgen\OKzero$, and let $\Ng$ be a non-zero integral ideal
of~$\OKzero$.
Suppose that there is a principally polarised ideal class~$T_1$
for~$(\Oc, \Phi)$.
By Proposition~\ref {prop:Acoprime} we may
assume that it is represented by a quadratic polynomial
$Q_1 = A_1 X^2+B_1 X+C_1 \in \OKzero [X]$ that is semiprimitive
modulo~$F \Ng$ with $\gcd(A_1,F \Ng)=1$.
Then there is an $\Ng$-system $\Qc = \{ Q_1, \ldots, Q_h \}$
for $(\Oc, \Phi)$ containing the given~$Q_1$.
\end{theorem}

\begin{proof}
Start with an arbitrary set of polynomials $\{ Q_1, \ldots, Q_h \}$
representing an orbit of the principally polarised ideal classes under
$\Cg_{\Oc, \Phi} (1)$.
Let $Q = A X^2 + B X + C$ with root $z$ be one of the $Q_i$ for $i \geq 2$.
Using Proposition~\ref {prop:Acoprime}, we may change $Q$ and~$z$ such that
they represent the same class, $Q$ is still semiprimitive modulo~$F \Ng$
and $\gcd (A, F \Ng) = 1$.

In the next step, we scale $Q$ to make it equiprimitive with~$Q_1$
modulo~$F \Ng$ while preserving the conditions on~$A$; the following lemma
(with $\mg = F \Ng$) provides the required scaling factor~$\epsilon$.
\begin {lemma}
\label{lem:equiprimitive}
Let $Q_1 = A_1 X^2 + B_1 X + C_1$ and
$Q = Q_2 = A_2 X^2 + B_2 X + C_2 \in \OKzero [X]$
be semiprimitive quadratic polynomials modulo some ideal~$\mg$ of~$\OKzero$,
representing principally polarised ideal classes for the same~$(\Oc, \Phi)$.
Then there is a (unique) $\epsilon \in K_0^\times$ such that
$Q_1$ and $A_2' X^2 + B_2' X + C_2' = \epsilon (A_2X^2+B_2X+C_2)$
are equiprimitive modulo~$\mg$.

Moreover, if $A_2$ is coprime to $\mg$, then so is~$A_2'$.
\end {lemma}

\begin {proof}[Proof of Lemma~\ref {lem:equiprimitive}]
With the notation of Lemma~\ref {lem:sqrtd}, let
$\epsilon = \delta_1 \delta_2^{-1}$; then the second polynomial, scaled
by $\epsilon$, has the same discriminant as the first polynomial.
But $\epsilon$ is in general not an algebraic integer, so it is a priori
not clear that the scaled polynomial still has integral coefficients.
However, we have $\dg_2' = \gcd (A_2', B_2', C_2') = \epsilon \dg_2 = \dg_1$,
which is an integral ideal, so $A_2'$, $B_2'$ and $C_2'$ are all integral.
Since $\epsilon$ is totally positive and coprime to~$\mg$ by
Lemma~\ref {lem:sqrtd}, semiprimitivity is preserved, and $A_2'$ remains
coprime to $\mg$ if $A_2$ is.
Unicity of $\epsilon$ is clear.
\end {proof}

Now \point{\ref{item:nsystem2}} of Definition~\ref {def:nsystem} is satisfied,
and we look for $\beta\in\OKzero$
such that $z' = z+\beta$ satisfies~\point{\ref{item:nsystem4}}.
Note that we then have $A' = A$, $B' = B-2A\beta$
and $D' = D$,
so the system remains equiprimitive and
\point{\ref{item:nsystem2}} remains satisfied.

Since by equiprimitivity the discriminants satisfy
$B^2 - 4 A C = B_1^2 - 4 A_1 C_1$,
we have $4 \mid B^2 - B_1^2 = (B - B_1)(B + B_1)$. From
$B - B_1 \equiv B + B_1 \pmod 2$ we deduce
$2 \mid B - B_1$.
For~\point {\ref{item:nsystem4}},
it suffices to take $\beta$
such that $A \beta \equiv \frac {B - B_1}{2} \pmod \Ng$,
which is possible since $\gcd (A, \Ng) = 1$.
This finishes the proof of Theorem~\ref {th:nsystems}.
\end{proof}

For the reader's convenience, we summarise the constructive proof of
Theorem~\ref {th:nsystems} in the following algorithm, before discussing
in more detail how single steps of it can be carried out.

\begin {algorithm}
\inoutput {A quadratic polynomial $Q_1$ as in Theorem~\ref {th:nsystems}}
{An $\Ng$-system $\Qc$ containing $Q_1$}
\label {alg:nsystem}
\begin {enumerate}
\item
Enumerate the orbit of the polarised ideal class represented by $Q_1$
under the action of $\Cg_{\Oc, \Phi} (1)$, and let $Q_1, \ldots, Q_h$ be the
resulting polynomials in $\OKzero [X]$.
\item
For $i = 2, \ldots, h$, write $Q_i = A_i X^2 + B_i X + C_i$
and do the following:
\begin {enumerate}
\item
Multiply $A_i$, $B_i$ and $C_i$ by an element in~$K_0$ such that they
become elements of~$\OKzero$ with $A_i \gg 0$ and
$\gcd (A_i, B_i, C_i, F \Ng) = 1$ as in
Proposition~\ref {prop:quadformexists}.
\item
Modify $Q_i$ by a matrix in~$\Sl_2 (\OKzero)$ as in
Proposition~\ref {prop:Acoprime} such that the new~$Q_i$
satisfies $\gcd (A_i, F \Ng) = 1$, while remaining semiprimitive
modulo~$F \Ng$.
\item
As in Lemma~\ref {lem:equiprimitive}, multiply $Q_i$ by
$\epsilon = \delta_1 \delta_i^{-1}$ with $\delta_i = 2 A_i z_i + B_i$.
\item
Let $\beta \in \OKzero$
be such that $A_i \beta \equiv \frac {B_i - B_1}{2} \pmod \Ng$;
replace $C_i$ by $A_i \beta^2 - B_i \beta + C_i$
and $B_i$ by $B_i - 2 \beta A_i$.
\end {enumerate}
\end {enumerate}
\end {algorithm}

The details of Step~\point {1} of the algorithm are out of the scope of
this article; see, for instance, \cite{EnTh14,Asuncion-thesis} for the
computation of the CM class group and the orbits.
An implementation is provided by \texttt {cmh}~\cite {cmh-1.1.1}, to which we
have added the remaining steps of
Algorithm~\ref {alg:nsystem} as the function \verb!nsystem!.

Step~\point {2a} is an application of strong approximation as in
\cite [Corollary~1.2.9]{Cohen00}:
Given a finite set $\Pc_0$ of prime ideals of~$\OKzero$
(the primes dividing $\gcd (A_i, B_i, C_i, F \Ng))$,
integers $e_\pg$ for $\pg \in \Pc_0$ (the negatives of the exponents of
$\pg$ in the $\gcd$),
$\Pc_\infty$ the set of the real embeddings of~$K_0$
and signs $e_v \in \{ \pm 1 \}$ for $v \in \Pc_\infty$
(the signs of $A_i$ under~$v$), we need to find
an element $\alpha \in \OKzero$ such that $v_\pg (\alpha) = e_\pg$
for $\pg \in \Pc_0$ and $\sign (v (\alpha)) = e_v$ for $v \in \Pc_\infty$.
In the PARI/GP system for number theory \cite {parigp}, for instance,
this can be implemented using the function \texttt {idealchinese}.

One way to do Step~\point {2b} is to construct matrices $\Mo_\pg$ as
in~\eqref {eq:Aprimetop} for every $\pg \mid F \Ng$.
Chinese remaindering (using again \texttt {idealchinese})
provides a matrix $\mmod {\Mo}{\mg} \in \Sl_2 (\OKzero / \mg)$, where
$\mg = \rad (F \Ng)$, with $\mmod {\Mo}{\mg} \equiv \Mo_\pg \pmod \pg$.
The question is now how to lift this matrix to $\Sl_2 (\OKzero)$.
We may start with an arbitrary lift
$\begin {psmallmatrix} a & b \\ c & d \end {psmallmatrix}$
to $\Mat_2 (\OKzero) \cap \Gl_2 (K_0)$.
If $a$ and $b$ are not coprime, then we may replace~$b$, using the Chinese
remainder theorem, by a $b' \in \OKzero$ satisfying $b' \equiv b \pmod \mg$
and $b' \equiv 1 \pmod {\gcd (a, b)}$. This is possible since
$\mg$ and $\gcd (a, b)$ are coprime: Otherwise, the determinant of
the matrix would have a non-trivial greatest common divisor with~$\mg$.
Then we compute
$u$ and $v$ in $\OKzero$ such that
$a u + v b = 1$ as follows:
By the Chinese remainder theorem, let $t \in \OKzero$ such that
$t \in a \OKzero$ and $t - 1 \in b \OKzero$, and let
$u = t / a$ and $v = (1 - t) / b$.
Write $m = a d - b c - 1 \in \mg$.
Then
$$\Mo =
\begin {pmatrix}
a       & b \\
c + v m & d - u m
\end {pmatrix}$$
has determinant~$1$ and reduction $\mmod {\Mo}{\mg}$ modulo~$\mg$.
This process is deterministic, polynomial time if the factorisation of
$F \Ng$ is known, and produces matrices~$\Mo$ with polynomial size
coefficients.
Alternatively, one may draw random matrices $\Mo \in \Sl_2 (\OKzero)$
until the resulting~$Q_i$ satisfies the desired coprimality conditions.

The value $\epsilon$ of Step~\point {2c} is most conveniently computed
as the totally positive square root of
$(B_1^2 - 4 A_1 C_1)(B_i^2 - 4 A_i C_i)^{-1}$.

Clearly Step~\point {2d} amounts to yet another application of the Chinese
remainder theorem.

\section {Complex conjugation}
\label{sec:cconj}

Class invariants $f (\tau)$ with $f \in \Fc_N$ are roots of class
polynomials of degree $h = |\Cg_{\Oc, \Phi} (1)|$ with coefficients
that lie \textit {a priori} in~$\Kr$. But it is well-known that
the $j$-invariant for $g = 1$ and the Igusa invariants for $g = 2$,
both of level $N = 1$, lead to class polynomials with coefficients in
the real subfield~$\Krzero$ of~$K$.
In this section we examine criteria under which this happens for higher
level modular functions and arbitrary~$g$ in our framework. This is
not only of theoretical interest, but also leads to a considerable speed-up
of algorithms computing class polynomials by floating point approximations.

As a preparatory step, we look at how complex conjugation acts on our
different algebraic structures and complex values.

\begin {lemma}
\label {lem:complexconjugation}
Assume the familiar setting of this article, that is,
$(\bg, \xi)$ and $z$, a root of the quadratic polynomial
$A X^2 + B X + C$,  are as in Propositions~\ref {prop:twogenerators}
and~\ref {prop:quadformexists}, the period matrix $\tau$ is computed
as in Proposition~\ref {prop:symplecticbasisgeneral}, and $f$
is a quotient of two Siegel modular forms with real $q$-expansions.

If $\Ac = \CC^g / \Phi (\bg)$ is the abelian variety associated to
$(\bg, \xi)$, then its complex conjugate variety $\cconj \Ac$ is induced
by $(\bg', \xi')$, where $\bg' = \cconj \bg = z' \OKzero + \OKzero$
with $z' = - \cconj z$ and $\xi' = \xi$.
The value $z'$ is a root of the quadratic polynomial
$A' X^2 + B' X + C'$ with $A' = A$, $C' = C$ and $B' = -B$,
and the associated period matrix is $\tau' = - \cconj \tau$.
Finally $f (\tau') = \cconj {f (\tau)}$.
\end {lemma}

\begin {proof}
The assertion on $\cconj \Ac$ follows from the fact that complex
conjugation commutes with the embeddings forming the CM type~$\Phi$,
see~\cite[Proposition~3.5.5]{langCM}.
The values of $\bg'$, $z'$, $A'$, $B'$ and $C'$ are clear by definition,
that of $\tau'$ is computed by
Proposition~\ref {prop:symplecticbasisgeneral}.

We now consider $q$-expansions.
To $\tau = \begin {psmallmatrix} \tau_1 & \tau_2 \\ \tau_2 & \tau_3 \end
{psmallmatrix}$
we associate the values $q_k = e^{2 \pi i \tau_k}$,
and similarly the $q_k'$ are associated to $\tau'$.
Then $q_k' = \cconj{q_k}$, and the $q$-expansions of the numerator and the
denominator of~$f$ having real coefficients implies that
$f (\tau') = \cconj {f (\tau)}$.
\end {proof}

In the following two sections we transpose two theorems on class
polynomials being defined over the real subfield from $g = 1$ to higher
dimension. To do so, we explicitly identify pairs of quadratic polynomials
that yield complex conjugate values. The main added difficulty is to examine
under which conditions these pairs belong to the same orbit under
$\Cg_{\Oc, \Phi} (1)$: If they belonged to different orbits, we would
obtain not one real class polynomial, but two complex conjugate class
polynomials, one for each orbit. We derive two separate, but related
criteria for the favourable situation
in Propositions~\ref {prop:cchyp} and~\ref {prop:cchyp2}.

\subsection {Real class polynomials for ramified levels}
\label {ssec:realramified}

The following result is the analogue of \cite [Theorems~4.4 and~6.1]{EnMo14};
it works for any function~$f$, but imposes severe restrictions on~$N$.

\begin {theorem}
\label {th:realramified}
Under the conditions of Theorem~\ref {th:conjugates},
assume furthermore that $F$ and $N$ are coprime and
that all primes dividing $N \OKzero$ are ramified in~$\OK$.
Suppose that the following holds for
the ideal $\bg=\bg_1=z_1\OKzero + \OKzero$:
\begin {equation}
\label{eq:cchyp}
\begin {aligned}
& \text {There exist $\cg \in \I (F)$ and $\mu\in K^\times$} \\
& \text {such that } \mu^{-1} \NphirO (\cg)^{-1} \bg = \cconj \bg
\text { and }
\mu \cconj \mu \Norm (\cg)=1.
\end {aligned}
\end {equation}
Then the class polynomial of Theorem~\ref {th:conjugates}
is an element of $\Krzero [X]$.

More precisely, let $\cg$ be as in \eqref{eq:cchyp}, let
(as in the theorem) $[\ag_i] \in \Cg_{\Oc, \Phi}(1)$
be such that
$[\ag_i]\cdot [(\bg, \xi)] = [(\bg_i,\xi_i)]$,
and let $c(i)$ be the index such that
$[\ag_{c (i)}] = \left[ \ag_i^{-1}\cg \right]$.
Then
\begin{equation}
\label{eq:formulaforccfinal}
\cconj{f(\tau_i)} = f (\tau_{c (i)}).
\end{equation}
\end {theorem}

Remark that the condition $N\mid C_1$ of Theorem~\ref{th:conjugates}
and the ramification condition imply by Theorem~\ref {th:n},
(3)$\Rightarrow$(1), that the ideal $N \OKzero$ is square-free.

\begin {proof}
Notice that it is sufficient to show~\eqref {eq:formulaforccfinal}, which
implies that the roots of the class polynomial are permuted by complex
conjugation and thus the coefficients are fixed by complex conjugation,
so that the polynomial is indeed an element of~$\Krzero [X]$.

Let $\Qc = \{Q_1,\ldots, Q_h\}$ be the $n$-system, and let
$Q_i = A_i X^2 + B_i X + C_i$ be one of its elements.
We first prove that $B_i$ is divisible by~$N$.
Let $\pg$ be a prime divisor of~$N \OKzero$. As $F$ and $A_i$ are coprime to~$N$,
while $C_i$ is divisible by~$N$ and $\pg$ is ramified in $K/K_0$, Lemma~\ref {lem:dedekind}
implies that $\tilde U (X) = X (X + B_i)$ is a square modulo~$\pg$, which
is only possible if $\pg \mid B_i$. Since $N \OKzero$ is square-free,
this shows that $N \mid B_i$.

Let $Q_i' = A_i' X^2 + B_i' X + C_i'$ with $A_i' = A_i$, $C'_i = C_i$,
$B_i' = -B_i$, with root $z_i' = -\cconj{z_i}$ and associated period matrix $\tau_i'$.
Then by Lemma~\ref {lem:complexconjugation} the left hand side
of~\eqref {eq:formulaforccfinal} equals $f (\tau_i')$, which to simplify
we will call the value of~$f$ in~$Q_i'$.
The quadratic polynomial $Q_i'$ represents the polarised ideal class
\begin {eqnarray*}
\left[\left( \cconj \bg_i, \xi_i \right) \right]\!\!
& = & [\cconj \ag_i]  [(\cconj \bg, \xi)]
\text { since the CM class group action is compatible with} \\
&& \text {complex conjugation} \\
& = & [\cconj \ag_i \cg] [(\bg, \xi)]
\text { by \eqref {eq:cchyp}.}
\end {eqnarray*}

The right hand side of~\eqref {eq:formulaforccfinal} is the value of~$f$
in~$Q_{c (i)}$, which represents the same polarised ideal class
$\left[ \ag_i^{-1} \cg \right] [(\bg, \xi)]
= \left[ \cconj \ag_i \cg \right] [(\bg, \xi)]$,
where we have used that $\cconj \ag_i \ag_i$, as an element of
$S_{\Oc, \Phi} (1)$ with $\mu = \Norm (\ag_i)$ in~\eqref {eq:ray},
is the identity of the CM class group.

It is now readily verified that $Q_i'$ satisfies the conditions of the
$N$-system as imposed by~$Q_i$:
The two polynomials are equiprimitive and satisfy
$\gcd (A_i', N) = \gcd (A_i, N) = \OKzero$
as they have the same $A_i' = A_i$ and discriminant;
furthermore $B_i - B_i' = 2 B_i$ is divisible by~$2 N$.

So $\Qc' = \Qc \backslash \{ Q_{c (i)} \} \cup \{ Q_i' \}$ is an
$N$-system.
If $h = |\Cg_{\Oc, \Phi}(1)| > 1$, then the two $N$-systems share a common
element and by Theorem~\ref {th:conjugates} lead to the same values of~$f$,
the roots of the class polynomial; so the values of $f$ in $Q_{c (i)}$ and
$Q'$ are the same, that is, \eqref {eq:formulaforccfinal} holds.
If $h = 1$, the proof of Theorem~\ref {th:conjugates} shows that $f$ has the
same value in the two quadratic polynomials representing the same polarised
ideal class and satisfying the $N$-system congruences.
\end {proof}

This proof is constructive in the sense that it allows to immediately
identify pairs of elements of the $N$-system that yield complex conjugate
values, which almost halves the time needed to compute floating point
approximations of the values of~$f$.

We need to examine more closely the validity of~\eqref {eq:cchyp}.
It is easily shown to hold for one element of an orbit if and only
if it holds for all elements of the orbit.
In practice, we will mainly consider situations in which it is known
to hold for all principally polarised ideals.

\begin{proposition}
\label {prop:cchyp}
Let $\Oc$ be an order in a CM field~$K$ with $\Oc\supseteq \OKzero$,
and let $\Phi$ be a primitive CM type of~$K$.
Let $(\bg, \xi)$ be a principally polarised ideal of $(\Oc,\Phi)$
such that $\bg$ is coprime to~$F$.
Then~\eqref{eq:cchyp}
holds for~$\bg$ in all of the following cases:
\begin{enumerate}
   \item \label{it:ccg2} $g\leq 2$;
   \item\label{it:g3} $g=3$ and $K$ contains an
      imaginary-quadratic subfield;
   \item\label{it:cpq} $g=6$, $\zeta_5\in \Oc$, the field $K$ is Galois over $\Q$,
      and the CM type $\Phi$ is CPQ-compatible
      as in~\cite{Somoza-thesis}.
\end{enumerate}
\end{proposition}
\begin{proof}
Suppose that \eqref{it:ccg2} holds.
If $g=1$, then one verifies directly that $\cg = \bg/\overline{\bg}$ and $\mu=1$ suffice.
If $g=2$, then take $\ag = \bg\OK$, $\cg = N_{\Phi}(\ag)\in \I(F)$,
and $\mu = \Norm (\ag)^{-1}$.
Lemma~I.8.4 of \cite{Streng10} then gives
\begin{equation}\label{eq:typenormoftypenorm}
N_{\Phi^r}(\cg) = \mu^{-1}\ag\cconj{\ag}^{-1}.
\end{equation}
Via the isomorphism between the groups of fractional ideals coprime to $F$ of $\Oc$ and~$\OK$,
we get $\NphirO (\cg) = \mu^{-1}\bg\cconj{\bg}^{-1}$,
which is the first equality of~\eqref{eq:cchyp}.
As we have $\Norm \left( \Nphi (\ag) \right) = \Norm (\ag)^g$
and thus $\Norm(\cg) = \Norm(\ag)^2$, we get $\mu\cconj{\mu}N(\cg) = 1$,
which is the second equality of~\eqref{eq:cchyp}.

In case \eqref{it:g3}, we can embed $K$ in such a way
into $\CC$ that $K = K^r$ (see e.g.~\cite[Propositions 3.3.2 and~3.4.2]{Kilicer-thesis}).
Then take $\ag = \bg\OK$, $\cg = \ag^{-1} \Nphi (\ag)\in \I(F)$,
and $\mu = \Norm(\ag)^{-1}$.
If $K$ is not Galois over~$\Q$, then \eqref{eq:typenormoftypenorm}
is exactly equation (3.4.2) of~\cite{Kilicer-thesis},
which can be shown to also hold in the Galois case using the
same argument.
The proof is finished exactly as for $g=2$.
Alternatively in the Galois case the displayed equation in the proof
of \cite[Proposition 3.3.5]{Kilicer-thesis}
gives a $\cg$ that works with $\mu=1$.

In case \eqref{it:cpq}, there are two possibilities for the
Galois group $G = \Gal (K/\Q)$ by \cite[Proposition~4.3.4]{Somoza-thesis}.

The first case is that $G$ is cyclic and
generated by $s$ of order~$12$, and after choosing an appropriate
embedding of $K$ into $\CC$ we have $K^r=K$,
$\Phi^r = \{1,s^4,s^7,s^8,s^9,s^{11}\}$,
the fixed field of~$s^4$ is~$\Q (\zeta_5)$ and
complex conjugation is $s^6$
by \cite[Proposition 4.3.11]{Somoza-thesis}.
Take $\ag = \bg \OK$ and $\alpha\in\Q(\zeta_5)$ such that
$\alpha\Z[\zeta_5] = N_{K/\Q(\zeta_5)}(s (\ag))$
and let $\mu = \cconj {\alpha} / \alpha \in\Q(\zeta_5)$.
Let $\cg = (s(\ag) s^5(\ag))/(\ag s^6(\ag))$.
A direct computation (or the first displayed equation of
the proof of \cite[Proposition 4.3.13]{Somoza-thesis})
again gives~\eqref{eq:typenormoftypenorm}.
We also have $\mu\cconj{\mu} = 1$ and $N(\cg)=1$.

The final case is where
$G = \left\langle s, t : ts=s^5t, t^2=s^3, s^6=1 \right\rangle$,
the fixed field of $s^2$ is $\Q(\zeta_5)$, complex conjugation is $s^3$,
and after choosing an appropriate
embedding of $K$ into $\CC$ we have $K^r=K$
and $\Phi^r = \{1,s^2,s^4,s^3t,s^4t,s^5t\}$
by \cite[Proposition 4.3.15]{Somoza-thesis}.
Take $\ag = \bg \OK$ and $\alpha\in\Q(\zeta_5)$ such that
$\alpha\Z[\zeta_5] = N_{K/\Q(\zeta_5)}(t(\ag))$
and let $\mu = \overline{\alpha}/\alpha$.
Let $\cg = t(\ag)/s^5t(\ag)$.
A direct computation (or the displayed equation of
the proof of \cite[Proposition 4.3.16]{Somoza-thesis})
again gives~\eqref{eq:typenormoftypenorm}.
We also have $\mu\cconj{\mu} = 1$ and $N(\cg)=1$.
\end{proof}

Notice that by the discussion following~\eqref {eq:normb} the fractional
ideal~$\bg_1$ occurring in Theorem~\ref {th:realramified} is coprime to~$F$
since it comes from a quadratic polynomial with $\gcd (A_1, F) = 1$,
so Proposition~\ref {prop:cchyp} applies in this case.

\subsection {Real class polynomials from the Fricke involution}
\label{ssec:fricke}

The following result is a generalisation of \cite[Theorem~3.4]{EnSc04};
it makes stronger assumptions on the function than
Theorem~\ref {th:realramified}, but does not require the primes dividing~$N$
to ramify. Again, it provides an explicit criterion for pairing up elements
of the $N$-system leading to complex conjugate roots of the class polynomial.

\begin {theorem}
\label{th:involutionS}
Under the conditions of Theorem~\ref {th:conjugates}, suppose furthermore
that we have $\gcd(N, C_i/N) = 1$ for all elements
$Q_i = A_i X^2 + B_i X + C_i$ of the $N$-system
and that $f$ is invariant under the Fricke involution
$\iota : \tau \mapsto -N \tau^{-1}$ of~$\Hc_g$.
Denote by $z_i$ the roots of the $N$-system and by $[(\bg_i, \xi_i)]$ the
associated polarised ideal classes. Let $\NNg = N \Oc + A_1 z_1 \Oc$.
Assume that the following hypothesis is satisfied for $\bg = \bg_1=z_1\OKzero+\OKzero$:
\begin {equation}
\label{eq:cchyp2}
\begin {aligned}
& \text {There exist $\cg \in \I (F)$ and $\mu\in K^\times$} \\
& \text {such that } \mu^{-1} \NphirO (\cg)^{-1} \bg \NNg = \cconj \bg
\text { and }
\mu \cconj \mu \Norm (\cg) = N.
\end {aligned}
\end {equation}
Then the class polynomial of $f (\tau_1)$ is an element of $\Krzero [X]$.

More precisely, let $[\ag_i] \in \Cg_{\Oc, \Phi}(1)$ be such that
$[\ag_i]\cdot [(\bg, \xi)] = [(\bg_i,\xi_i)]$,
and let $c (i)$ be the index such that
$[\ag_{c (i)}] = \left[ \ag_i^{-1}\cg \right]$.
Then
\begin{equation}
\label{eq:formulaforcc2final}
\cconj{f(\tau_i)} = f (\tau_{c (i)}).
\end{equation}
\end{theorem}

Before proving this theorem, we have a closer look at the different
assumptions occurring in it, assessing how to verify them and whether
they represent important limitations.

The congruence condition $\gcd(N, C_i/N) = 1$ is not essential:
By Theorem~\ref {th:n}\point {2} we may assume that it holds for $i=1$,
and then impose that $B_i = B_1$ modulo higher powers of prime ideals
dividing $N \OKzero$ (very crudely, switching to an $N^2$-system will work).

Condition~\eqref{eq:cchyp2} is a variation on~\eqref{eq:cchyp}.
In fact,
letting $N = 1$ and thus $\NNg = \Oc$ in~\eqref {eq:cchyp2}
retrieves~\eqref {eq:cchyp} as a special case.
Moreover, if~\eqref {eq:cchyp} is valid (for instance under the conditions
of Proposition~\ref {prop:cchyp}), then~\eqref {eq:cchyp2} becomes
equivalent to the following:
\begin {equation}
\label{eq:Nhyp}
\begin {aligned}
& \text {There exist $\cg \in \I (F)$ and $\mu\in K^\times$} \\
& \text {such that } \NNg = \mu \NphirO (\cg)
\text { and }
\mu \cconj \mu \Norm (\cg) = N,
\end {aligned}
\end {equation}
which again shows that~\eqref {eq:cchyp2} is independent of the element
in the orbit, at least under~\eqref {eq:cchyp}, and that in reality it is
a statement about~$\NNg$.
We will see in Proposition~\ref{prop:cchyp2} that this condition often holds.
The example in Section~\ref{ssec:cchyp2} shows that it
is necessary in the sense that the theorem becomes false
when replacing the condition with~\eqref{eq:cchyp}.

In order to prove the theorem, we have a look at several involutions
occurring in our context, such as the Fricke involution $\iota$ and
the involutions $z\mapsto -N/z$ and $ z \mapsto -1 / z$.

\begin{lemma}
\label{lem:S}
Let $\Z$-bases $\Bc_1$ and $\Bc_2$ of~$\OKzero$ be given as in
Proposition~\ref {prop:symplecticbasisgeneral}, and let
$S \in \Gl_g (\Z)$ be the matrix such that
$\T{\Bc_2} = S\T{\Bc_1}$.
Then $S$ is symmetric.

If $\tau$ is the period matrix associated to $z$ by
Proposition~\ref {prop:symplecticbasisgeneral},
then the period matrix associated to $-1/z$ is
\begin{equation}
\label{eq:inverse}
\tau' =
\begin{pmatrix}
0 & -S^{-1} \\ S & \phantom{-}0^{\phantom{-1}}
\end{pmatrix}
\tau =
- (S \tau S)^{-1}.
\end{equation}
\end{lemma}

\begin{proof}
Notice that $S = [1]_1^2$ and $S^{-1} = [1]_2^1$ with the notation
of~\eqref {eq:blowup}. Then by Lemma~\ref{lem:M}, the matrix
$$
\Mg = \otogleft{\begin{pmatrix} 0 & -1 \\ 1 & 0\end{pmatrix}}
= \begin{pmatrix} 0 & -S^{-1} \\ S & \phantom{-}0^{\phantom{-1}}\end{pmatrix}
$$
satisfies $\Mg\in \Sp(\Z)$, which implies $S=\T {S}$.

It remains to show that
$\tau' = \Mg \tau$,
which is the special case
$M = \begin {psmallmatrix} 0 & -1 \\ 1 & 0 \end {psmallmatrix}$
of Lemma~\ref{lem:MzMtau}.
\end{proof}

\begin {example}
For $g = 2$ and $\diffgen = \sqrt {\Delta_0}$, let $\Bc_1$ and $\Bc_2$
be as in the proof of Corollary~\ref {cor:symplecticbasisg2}.
Then
$$S = \begin {pmatrix} 0 & -1 \\ -1 & 0 \end {pmatrix}\
\mbox{if $\Delta_0$ is even, and}\
S = \begin {pmatrix} 0 & -1 \\ -1 & 1 \end {pmatrix}\
\mbox{if $\Delta_0$ is odd.}$$
\end {example}

\begin{lemma}\label{lem:involutions}
There is an involution $\iota':\Hc_g\rightarrow \Hc_g$
with the following properties.
\begin{enumerate}
   \item If $\tau$ and $\tau'$ correspond to
      respectively $z$ and $-N / z$,
      then $\tau' = \iota'(\tau)$.
   \item
      Let $f$ be a modular function for $\Gamma^0(N)$.
      Then $f$ is invariant under the Fricke involution
      $\iota:\tau\mapsto -N\tau^{-1}$
      if and only if it is invariant under $\iota'$.
   \item
      If $g=1$, then $\iota=\iota'$.
\end{enumerate}
\end{lemma}
\begin{proof}
Let $S = [1]_1^2$ be as in Lemma~\ref{lem:S}.
Let $\iota' = \iota_S : \Hc_g\rightarrow \Hc_g$ be the involution
given by $\tau \mapsto -N (S\tau S)^{-1}$.
Note that for $g=1$ we have $\Bc_1=\Bc_2 = (1)$, hence $S=1$
and $\iota'=\iota$, which proves the third statement.

The first statement is \eqref{eq:inverse} in Lemma~\ref{lem:S}.

Finally,
we have $\iota^2=\mathrm{id}_{\Hc_g}$
and hence
$$\iota^{-1} \iota' \tau = \iota\iota'\tau = -N(-N(S\tau S)^{-1})^{-1}
= S\tau S =\begin {pmatrix} S & 0 \\ 0 & S^{-1} \end {pmatrix} \tau.$$
As the latter matrix is in $\Gamma^0(N)$, we find that $f$ is
invariant under $\iota$ if and only if it is
invariant under $\iota'$.
\end{proof}

\begin {proof}[Proof of Theorem~\ref{th:involutionS}]
We follow the same approach as in the proof of
Theorem~\ref {th:realramified}.

By Lemma~\ref{lem:complexconjugation}, we have
$\cconj {f(\tau_i)} = f(-\cconj {\tau_i})$, with the latter period
matrix corresponding to $- \cconj z_i$.
By Lemma~\ref{lem:involutions}\point {2} and invariance of~$f$ under
the Fricke involution,
$f (-\cconj {\tau_i}) = f (\iota'(-\cconj {\tau_i}))$,
and by~\ref{lem:involutions}\point {1} the latter period matrix
corresponds to $z_i' = N / \cconj z_i$.
This is a root of the polynomial $Q_i' = A_i' X^2 + B_i' X + C_i'$
with $A_i' = C_i / N$, $B_i' = B_i$ and $C_i' = A_i N$,
which satisfies the $N$-system conditions as imposed by~$Q_i$:
The discriminant $B_i^2 - 4 A_i C_i$ of~$Q_i$ is totally negative, so that
with $A_i$ totally positive also $C_i$ is totally positive. Since
$\gcd(N, C_i/N) = 1$ by assumption, the polynomial~$Q_i'$ is semiprimitive
modulo~$N$ as in Definition~\ref {def:semiprimitive}.
Moreover, the discriminants of $Q_i'$ and $Q_i$ are equal,
which implies equiprimitivity modulo~$N$ as in
Definition~\ref{def:equiprimitive}.

The left hand side of~\eqref {eq:formulaforcc2final} is the value of~$f$
in~$Q_i'$, the right hand side is the value of~$f$ in $Q_{c (i)}$. If we
can show that the two polynomials represent the same principally
polarised ideal class, we conclude as in the proof of
Theorem~\ref {th:realramified} that they lead to the same value and thus
finish the proof of Theorem~\ref {th:involutionS}.

Condition \eqref{eq:cchyp2} and invertibility of~$\bg$
(Proposition~\ref{prop:quadformexists}) imply that $\NNg$
is invertible and
$[\cg] \cdot [(\bg_1, \xi_1)]
= \left[ (\NNg^{-1} \cconj{\bg_1}, N \xi_1) \right]$,
and this gives us a handle on~$Q_{c (i)}$, which represents
the polarised ideal class
\begin {eqnarray*}
[\ag_{c (i)}] [(\bg_1, \xi_1)]
& = & \left[ \ag_i^{-1} \cg \right] [(\bg_1, \xi_1)]
   \text { by definition} \\
& = & \left[ \cconj \ag_i \right]
      \left[ (\NNg^{-1} \cconj \bg_1, N \xi_1) \right] \\
&& \text {from the action of $\cg$ and since
      $\left[ \ag_i^{-1} \right] = \left[ \cconj \ag_i \right]$
      as seen before} \\
& = & \left[ (\NNg^{-1} \cconj \bg_i, N \xi_i) \right]\\
&& \text {by compatibility of the type norm with complex conjugation.}
\end {eqnarray*}
It remains to show that $Q_i'$ represents the same polarised ideal class.
By Proposition~\ref {prop:twogenerators} it leads to the
polarised ideal
$(\bg_i', \xi_i')$ with
$\bg_i' = (N / \cconj z_i) \OKzero + \OKzero$ and
\[
\xi_i'
= \left( (N / \cconj{z_i} - N / z_i) \lambda \right)^{-1}
= N^{-1} z_i \cconj z_i \left( (z_i - \cconj z_i) \lambda \right)^{-1},
\]
which is equivalent to
$\left( N^{-1} \cconj z_i \bg_i', N^2/(z_i\cconj{z_i}) \xi_i' \right)=
\left( N^{-1} \cconj z_i \bg_i', N \xi_i \right)$.

The component $N\xi_i$ is already what we are looking for, and it remains to
show that
\begin {equation}
\label {eq:compfricke}
\NNg (N^{-1} \cconj z_i  \bg_i')
= \cconj \bg_i
\end {equation}
on the ideal side. We first examine $\NNg$ more closely.
Let $\theta_1 = A_1 z_1$ and $\theta_i = A_i z_i$. Since by equiprimitivity
the quadratic polynomials $Q_1$ and $Q_i$ have the same discriminant,
we obtain $\theta_1 - \theta_i = \frac {B_i - B_1}{2}$, and this difference
is divisible by~$N$ from the $N$-system congruences. It follows that
$\NNg = N \Oc + \theta_1 \Oc = N \Oc + \theta_i \Oc$.
Using
\[
\Oc = \OKzero + \dg_i^{-1} \theta_i
\]
of~\eqref {eq:orderfromform} we compute
\begin {eqnarray}
\NNg
& = & N \OKzero + N \dg_i^{-1} \theta_i
      + \OKzero \theta_i + \dg_i^{-1} \theta_i^2 \nonumber\\
& = & N \OKzero + N \dg_i^{-1} \theta_i + \OKzero \theta_i
      + \dg_i^{-1} (B_i \theta_i + A_i C_i) \nonumber\\
& = & N \OKzero + N \dg_i^{-1} \theta_i + \OKzero \theta_i
      \text { since } N | C_i ,\ \dg_i | A_i \text { and } \dg_i | B_i \nonumber\\
& = & N \OKzero + \dg_i^{-1} \theta_i
      \text { since } \gcd (\dg_i, N) = 1.\label{eq:formulaN}
\end {eqnarray}
So the left hand side of~\eqref {eq:compfricke} is computed as
\begin{align*}
&  (N \OKzero + \dg_i^{-1} A_i z_i)(\OKzero + N^{-1} \cconj z_i \OKzero) \\
&= N \OKzero + \cconj z_i \OKzero + \dg_i^{-1} A_iz_i + \dg_i^{-1} C_i / N \\
&= N \OKzero + \cconj z_i \OKzero + \dg_i^{-1} (B_i - A_i \cconj z_i) + \dg_i^{-1} C_i / N \\
&= N\OKzero + \cconj z_i \OKzero + \dg_i^{-1} B_i + \dg_i^{-1} C_i / N \\
&= \OKzero + \cconj z_i \OKzero \\
&= \cconj \bg_i,
\end{align*}
which finishes the proof.
\end {proof}

\begin {proposition}
\label {prop:cchyp2}
Let $K$ be a CM field with a primitive CM type $\Phi$. Let $z\in K\setminus K_0$
be a root of a quadratic polynomial $Q = AX^2 + BX + C\in\OKzero[X]$
and let $\Oc$ be the ring of multipliers of $\bg = z\OKzero + \OKzero$.
Suppose that $N\mid C$ and that all of $A$, $\gcd(B,C/N)$, and the
conductor $\fg$ of~$\Oc$ are coprime to~$N$.
In each of the following cases, the ideal $\NNg = \Oc + Az \Oc$ satisfies
\eqref{eq:Nhyp} and \eqref {eq:cchyp2}:
\begin{enumerate}
\item $g=1$;
\item $g=2$ and $N$ is coprime to the discriminant of $K$.
\end{enumerate}
\end{proposition}

We have identified more cases in which~\eqref {eq:Nhyp} holds, with somewhat
lengthy proofs relying only on Galois theory and unrelated to class
invariants; to lighten the presentation, we delegate these results
to another venue.

\begin {proof}
Since $g \leq 2$ the conclusion of Proposition~\ref {prop:cchyp} holds and
\eqref {eq:cchyp2} and \eqref {eq:Nhyp} are equivalent; we proceed to prove
the latter.

Letting $\theta = A z$ and $\dg = \gcd (A, B, C)$, we obtain by
\eqref {eq:formulaN} that
\begin{align}
\NNg \cconj{\NNg}
& = (N \OKzero + \dg^{-1} \theta) (N \OKzero + \dg^{-1} \cconj {\theta})
   \nonumber\\
& = N \left( \dg^{-1} \theta + N \OKzero + \dg^{-1} B
    + (\dg^{-1} A) (\dg^{-1} C/N) \right) \nonumber\\
& = N (\dg^{-1} \theta + \OKzero)
    \text { because of the conditions of coprimality with $N$}
\nonumber \\
&= N\Oc.
\label{eq:NNgNNgbarisNg}
\end{align}

If $g=1$, then without loss of generality
$\Nphir$ is the identity on $K = \Kr$;
we take $\mu = 1$ and $\cg = \NNg\OK$ and conclude
by~\eqref {eq:NNgNNgbarisNg}.

If $g=2$ and $N$ is coprime to the discriminant of~$K$, then
all primes $p\mid N$ are unramified in~$K$.
This and the assumption that $N$ is coprime to~$\fg$
imply that $B$ is coprime to~$N$ by Lemma~\ref{lem:dedekind}.
Then
\begin {equation}
\label {eq:Nsplit}
\NNg + \cconj{\NNg}
= N \OKzero + \dg^{-1} \theta + \dg^{-1} B
= \OKzero + \dg^{-1} \theta = \Oc.
\end {equation}
For every $p\mid N$, let $\NNg_p = \NNg + (p)$;
then $\NNg_p \cconj{\NNg_p} = (p)$,
and by the multiplicativity of~\eqref{eq:Nhyp}
it suffices to show the result for $p$ and $\NNg_p$
in the place of $N$ and $\NNg$.
We choose $\mu = 1$; then it is enough to show the existence of an
ideal~$\cg$ such that $\NphirO (\cg) = \NNg_p$, from which
$\Norm (\cg) = \NphirO (\cg) \cconj {\NphirO (\cg)}
= \NNg_p \cconj {\NNg_p} = p \OKzero$ follows.

Let $\pg$ be a prime ideal above $p$ of $L = K \Kr$, which is
the Galois closure of~$K$.
Theorems 1 and~2 of \cite{Goren97} give all possible decomposition groups
of~$\pg$ when $p$ is unramified in~$K$.
By~\eqref {eq:Nsplit} the primes
of~$\OKzero$ above~$p$ split in $K / K_0$, hence
only case (1) of Theorem~1 and cases (1), (3), and (5) of Theorem~2 occur.

Case (1) of Theorem~1 is the totally split case with
$\Gal(K/\Q)=\langle y\rangle$ for an element $y$ of order~$4$
such that $\Phi = \{1, y\}$, and $y^2$ is complex conjugation.
We have $(p) = \pg \cdot y(\pg) \cdot y^2(\pg) \cdot y^3(\pg)$
and $(p) = \NNg_p \cdot y^2(\NNg_p)$,
so that $\NNg_p = y^a (\pg) y^{a-1}(\pg)$ for
some~$a \in \{ 0, 1, 2, 3 \}$.
Without loss of generality, we may assume $a = 0$
(or, otherwise said, we may replace $\pg$ by $y^a (\pg)$);
then $\Phi^r = \{1, y^{-1}\}$ implies
$\Norm_{\Phi^r,\Oc}(\pg) = \NNg_p$.

Case (1) of Theorem~2 is also totally split, but now
$\Gal(L/\Q) = \langle x,y\rangle\cong D_4$
with $x^2 = y^4 = 1$ and $xyx = y^3$,
$\Gal(L/K) = \langle x\rangle$, $\Gal(L/K^r) = \langle xy^3\rangle$,
and $y^2$ is complex conjugation.
The eight prime ideals of $L$ above~$p$ are the $g(\pg)$ for $g\in D_4$
and the prime ideals of $K$ above~$p$ are the $y^a(\pg)\cdot xy^a(\pg)$
for $a=0,1,2,3$, so
$\NNg_p = y^a(\pg)\cdot xy^a(\pg)\cdot  y^{a-1}(\pg)\cdot xy^{a-1}(\pg)$
for some~$a$.
Without loss of generality (that is, after potentially choosing a
different~$\pg$), we can assume $a=0$, and then $\NNg_p$ is exactly what
appears as the type norm
$\NphirO (\pg\cap \Kr) = (\pg \cdot x(\pg))(y^3(\pg) xy^3(\pg))$ of Case~(1)
on page~38 of~\cite{Goren97}, just below Theorem~2.

Cases (3) and (5) of Theorem~2 are conjugate, so without loss of generality
(choosing a different $\pg$ for the same~$p$ if needed) we are in Case~(5).
In that case, pages 38 and~39 of~\cite{Goren97} give the decomposition
of $p$ in $K$ as a product of two primes
$(p) = (\pg y^3(\pg))(y(\pg) y^2(\pg))$,
and the type norm as
$\NphirO (\pg\cap \Kr) = (\pg y^3(\pg))$.
But then the prime $(\pg y^3(\pg))$ of $K$ is $\NNg_p$ or~$\cconj{\NNg_p}$, so that taking
$\cg = \pg\cap \Kr$ or $\cg=\cconj{\pg\cap \Kr}$ gives the desired result.
\end {proof}

\begin{remark}\label{rmk:motivationcchyp2}
We now sketch a more geometric view of Theorems
\ref{th:realramified} and~\ref{th:involutionS}
and their proofs,
giving an alternative explanation of where
the involution and the
complex conjugation hypotheses
\eqref{eq:cchyp} and~\eqref{eq:cchyp2} come from.

The variety $\Gamma^0(N)\backslash \Hc_g$ is the coarse moduli space of $N$-isogenies,
in the sense that its points correspond to triples
$t = (A, A', \varphi)$, where $A$ and $A'$ are principally polarised
abelian varieties of dimension~$g$
and $\varphi : A\rightarrow A'$ is an $N$-isogeny.
Each $z = z_i\in K\setminus K_0$ gives rise to such a triple $t = t_i$ as follows.
With $\bg = z\OKzero + \OKzero$ and $\xi$ as in Proposition~\ref{prop:twogenerators},
we get a principally polarized abelian variety $A/\overline{\Q}$ with $A(\CC) \cong \CC^g/\Phi(\bg)$,
which has CM by $(\Oc, \Phi)$ via
an isomorphism $\iota : \Oc\rightarrow \mathrm{End}(A)$.
Similarly from $z' = N^{-1} z$ we get $\bg'$, $\xi'$, $A'$, and~$\iota'$.
Moreover, we now have the $N$-isogeny $\varphi : A \rightarrow A'$ that is the
identity on~$\CC^g$.

Let $f$ be a function on $\Gamma^0(N)\backslash\mathbf{H}_g$ defined over~$\mathbf{Q}$,
as in our theorems.
The values $f(\tau_i) = f(t_i)$ will depend only on the isomorphism
classes of the corresponding triples~$t_i$.
Moreover, for any $\sigma\in\mathrm{Gal}(\overline{\Q}/\Q)$,
including complex conjugation,
we have $\sigma(f(t)) = f(\sigma(t))$.

To get a class polynomial with real coefficients,
we want to have $\overline{t} \cong \sigma(t)$
for some $\sigma\in \Gal(\overline{\Q}/K^r)$.
To see when this can and cannot happen we look
more closely at the kernel of~$\varphi$ as follows.

Under some mild assumptions,
we have (similarly to~\eqref{eq:compfricke}) the equality
$\bg' = \smash{\overline{\NNg}}^{-1}\bg$
for some ideal $\NNg$ with $\NNg\overline{\NNg} = (N)$.
It follows that $\varphi$ is actually an $\overline{\NNg}$-\emph{multiplication}
with respect to the CM type~$\Phi$,
that is, an isogeny satisfying
\begin{equation}\label{eq:defofNgbarmul}
\ker(\varphi) = A[\iota(\overline{\NNg})]:=
\bigcap\limits_{\beta\in\overline{\NNg}} \ker(\iota(\beta))
\text { and }
\iota'(\alpha)\circ \varphi = \varphi\circ\iota(\alpha)
\end{equation}
for all $\alpha\in\Oc$, where both $\iota$ and $\iota'$ have type~$\Phi$.

Complex conjugation turns $t=(A, A, \varphi)$ into $\overline{t} = (\overline{A}, \overline{A'},
\overline{\varphi})$.
The complex conjugate abelian varieties $\overline{A}$ and $\overline{A'}$
again have CM type $\Phi$ once we adorn them with
$\overline{\iota} : \alpha \mapsto \overline{\iota(\overline{\alpha})}$
and the similarly defined~$\overline{\iota'}$.
Then $\overline{\varphi}$ is an $\NNg$-multiplication
(with $\NNg$, not $\overline{\NNg}$)
as can be seen by taking complex conjugates of~\eqref{eq:defofNgbarmul}.
%
%

On the other hand, all $\Gal(\overline{\Q}/K^r)$-conjugates of
the $\overline{\NNg}$-multiplication~$\varphi$
are again $\overline{\NNg}$-multiplications with respect to the CM type~$\Phi$
(\cite[Proof of Prop.~11 in \S14.7;
see also \S7.6 and Prop.~31 of \S8.5]{shimura-taniyama}).
So in order for $\overline{t}$ to be $\Gal(\overline{\Q}/K^r)$-conjugate
to~$t$, we need to have $\NNg = \overline{\NNg}$.

This explains why Theorem~\ref{th:realramified}
assumes that all primes dividing $N$
ramify in~$K$: so that we have $\NNg = \overline{\NNg}$.
As the $\Gal(\overline{\Q}/K^r)$-action
is not always transitive, we need to supplement
this assumption with the assumption
that the first components $\overline{A}$
and $A$ of $\overline{t}$ and $t$ are in the same orbit,
which is exactly~\eqref{eq:cchyp}.

However, if $\NNg\not=\overline{\NNg}$,
then we need to somehow
turn our $\NNg$-multiplication $\overline{\varphi}$
back into an $\overline{\NNg}$-multiplication.
We can do so by taking its dual,
see \cite[\S14.4, Prop.~6]{shimura-taniyama},
which corresponds to taking the Fricke involution on the moduli space.
This is why we assume that $f$ is fixed under
the Fricke involution in Theorem~\ref{th:involutionS}:
so that $f$ makes no distinction between $\overline{t}$
and its dual
$(\overline{A'}, \overline{A}, \overline{\varphi}^\dagger)$.
To get a class polynomial with real coefficients,
we furthermore need $\overline{A'}$
and $A$ to be in the same
$\mathrm{Gal}(\overline{\Q}/K^r)$-orbit, which
is exactly~\eqref{eq:cchyp2} in Theorem~\ref{th:involutionS}.
\end{remark}

\section{Families of functions for
\texorpdfstring {$\Gamma^0 (N)$}{Gamma0 (N)}}
\label {sec:functions}

In this section we provide a few examples of families of functions
for $g=2$ that
can be used in the context of Theorem~\ref{th:maintheorem},
i.e., functions that are invariant under $\Gamma^0(N)$
and quotients of modular forms with rational $q$-expansions.
We will use them in Section~\ref{sec:examples} to provide numerical examples.
These are just a few examples. We expect that many more good functions exist,
and we leave a thorough search for future research.

\subsection{Functions obtained from Igusa invariants}

Igusa defines modular forms $h_4$, $h_6$, $h_{10}$ and $h_{12}$
with rational $q$-expansions that generate the graded ring of modular forms
for $\Spfour(\Z)$ \cite{igusa-i}; so for
$\Mg =
\begin {psmallmatrix} a & b \\ c & d \end {psmallmatrix}
\in \Spfour (\Z)$ and
$\Mg \tau = (a \tau + b)(c \tau + d)^{-1}$, one has
$h_k (\Mg \tau) = \det (c \tau + d)^k h_k (\tau)$.
Taking quotients of forms of the same
weight yields modular functions for $\Spfour (\Z)$ such as
$$i_1 = \frac{h_4h_6}{h_{10}},\quad
i_2=\frac{h_4^2h_{12}^{\phantom{2}}}{h_{10}^2},\quad\mbox{and}\quad
i_3 = \frac{h_4^5}{h_{10}^2}$$
known as \emph{absolute Igusa invariants}.
There are many possible choices of absolute Igusa invariants;
the above functions correspond to~\cite{Streng-runtime}.
Since these are modular functions for the full modular group,
their CM values are automatically class invariants.

Alternatively, one may take \emph{simple $h_k$-quotients}
$$\frac{h_k(\tau/N)}{h_k(\tau)}$$
stable under $\Gamma^0(N)$
or \emph{double $h_k$-quotients}
\begin{equation}
\label{eq:doublequot}
f =
\frac{h_k(\tau/N_1) h_k(\tau/N_2)}{h_k(\tau) h_k(\tau/(N_1 N_2))}
\end{equation}
stable under $\Gamma^0(N)$ for $N = N_1N_2$.
The latter function is also invariant under the Fricke
involution
$\iota : \tau \mapsto -N \tau^{-1}$ of
Theorem~\ref {th:involutionS}:
\[
f (\iota (\tau))
= \frac {h_k (- N_1 \tau^{-1}) h_k (- N_2 \tau^{-1})}
{h_k (- \tau^{-1}) h_k (- N_1 N_2 \tau^{-1})}
= f (\tau),
\]
where we have used
$h_k (- \tau^{-1}) = h_k (\Jg \tau) = \det(-\tau)^{k} h_k (\tau)$
for $\Jg$ as in~\eqref {eq:J}.

As for simple and double eta quotients in dimension~$1$, the process
may be generalised to obtain multiple quotients of~$h_k$,
cf.~\cite {EnSc13}.

For $k=10$, similar functions and their square roots
have also been studied by
de Shalit and Goren~\cite{deshalit-goren}.

\subsection{Theta products}
\label {ssec:thetaproducts}

The \emph{theta constant} of characteristic
$(\alpha,\beta)\in (\Q^g)^2$ is given by
$$\theta[\alpha, \beta](\tau)
= \sum_{n\in\Z^g} \exp \left( \pi i \T{(n+\alpha)} \tau (n+\alpha)
+ 2\pi i \T{(n+\alpha)} \beta \right)$$
for $\tau\in\Hc_g$.
For $\alpha$, $\beta\in\{0,1/2\}^g$, it is a modular form
of weight $\frac{1}{2}$
for $\Gamma(8)$ with $q$-expansion coefficients in~$\Q$.

From now on, we consider only the case $g=2$, and
we also use the abridged notation
\[
\theta_{8 a_1 + 4 a_2 + 2 b_1 + b_2}
= \theta \left[ \begin {pmatrix} a_1/2 \\ a_2/2 \end {pmatrix},
\begin {pmatrix} b_1/2 \\ b_2/2 \end {pmatrix} \right]
\]
introduced in \cite [\S6.2]{Dupont06}
for $a_1$, $a_2$, $b_1$, $b_2 \in \{ 0, 1 \}$.

Ibukiyama has shown in \cite [Theorem~A]{Ibukiyama91} that the graded ring
of modular forms for $\Gamma_0 (2)$ is generated by the four forms with
rational $q$-expansions given by
\begin {align*}
x &= (\theta_0^4 + \theta_1^4 + \theta_2^4 + \theta_3^4) / 4 \\
y &= (\theta_0 \theta_1 \theta_2 \theta_3)^2 \\
z &= (\theta_4^4 - \theta_6^4)^2 / 2^{14} \\
k &= (\theta_4 \theta_6 \theta_8 \theta_9 \theta_{12} \theta_{15})^2 / 2^{12}
\end {align*}
of respective weights $2$, $4$, $4$ and $6$;
notice that $2^{12} \, y k = h_{10}$.

Evaluating these forms in $\tau / 2$, we obtain generators for the graded
ring of modular forms for $\Gamma^0 (2)$ as
$X (\tau) = x (\tau/2)$, $Y (\tau) = y (\tau/2)$, $Z (\tau) = z (\tau/2)$
and $K (\tau) = k (\tau/2)$.
The smallest weight for which the vector space of forms
has dimension at least~$2$ is~$4$, with a basis given by $X^2$, $Y$ and $Z$.
By taking a quotient of two such forms, we obtain a function for
$\Gamma^0 (2)$, which we expect to yield small class invariants.
In fact, the second part of the theorem by Ibukiyama shows that the
field of Siegel modular functions for $\Gamma^0 (2)$ is rational of
transcendence degree~$3$ and generated by $X^2 / Y$, $Z / Y$ and $X^3 / K$.

We may also fix $F$ as one of $X$, $Y$, $Z$ or $K$ and consider simple
quotients
$\frac {F (\tau / N)}{F (\tau)}$, which are functions for
$\Gamma^0 (2 N)$, and double quotients
$\frac {F (\tau / N_1) F (\tau / N_2)}{F (\tau) F (\tau / (N_1 N_2))}$,
which are functions for $\Gamma^0 (2 N_1 N_2)$.
The form $F = X$ is
promising in this context due to its low weight,
as is the divisor $F=Y$ of~$h_{10}$.

\section{Numerical examples}
\label{sec:examples}

\subsection {Implementation and setup for reproducibility}

We have implemented the algorithms described above in the PARI/GP
system~\cite {parigp-2.15.4} at version 2.15.4; the code for reproducing
our example class polynomials is made available as supplementary material
to this article~\cite {nsystemcmh-0.2}.
More precisely, we have added code for computing $N$-systems to the GP
script \texttt {shimura.gp} in the latest release 1.1.1 of
CMH~\cite {cmh-1.1.1}, developed by the first author, Emmanuel Thomé
and Régis Dupont and described in~\cite {EnTh14}.
This $N$-system code is now also used in CMH for Igusa class
polynomials with $N=1$.
The CMH library provides C code for asymptotically fast evaluations of
Siegel modular forms using arbitrary precision floating point operations.
It can be called from GP scripts using the
PariTwine software~\cite {paritwine-0.2.0} at version 0.2.0,
which relies additionally on
GNU MP at version 6.2.1~\cite {gmp-6.2.1},
GNU MPFR at version 4.2.0~\cite {mpfr-4.2.0}
and
GNU MPC at version 1.3.1~\cite {mpc-1.3.1}.
The GP script \texttt {common.gp} modifies the $N$-system computation of
CMH to enforce the divisibility conditions of Theorems~\ref {th:maintheorem}
and~\ref {th:involutionS} of the values $C_i$ using the algorithms
of the proof of Theorem~\ref {th:n}. It then calls CMH through
PariTwine to compute the complex roots of the class polynomial and
it guesses the correct algebraic class polynomial.
It is called by the scripts \texttt {example1.gp}, \texttt {example2.gp}
and \texttt {example3.gp}, which compute the Igusa class polynomials and
our alternative class polynomials for the examples given in detail in
the following three subsections.

Additionally, we have used the second author's RECIP~\cite{recip-3.4.1}
code at version 3.4.1, which is developed as a package for
SageMath~\cite{sagemath-10.2}, to compute Igusa class polynomials;
its \texttt {class\_polynomial} command returns a result proved to be
correct by the approaches of Bouyer--Streng~\cite{BouyerStreng}
and Lauter--Viray~\cite{LauterViray}. In all cases, these Igusa polynomials
were identical to those computed using CMH and our GP scripts.

We followed the reproducibility approach described in~\cite {CoFeHiSw23}
to easily make the results available to the reader through
Guix~\cite {guix-bf17a01}. For reproducing the computation of the first
example, it is enough to run from the subdirectory \texttt {parigp}
of~\cite {nsystemcmh-0.2} the command
\begin {verbatim}
guix time-machine -C channels.scm -- shell -C -m manifest.scm -- \
   gp < example1.gp
\end{verbatim}
and analogously for the other examples, using the provided files
\texttt {channels.scm} and \texttt {manifest.scm}, which record
the computing environment with the exact dependencies of all
required software.

Outside of Guix, the reader may install the software packages given above
at their respective versions into the prefix \texttt {/usr/local}, say,
and then call
\begin {verbatim}
export GUIX_ENVIRONMENT=/usr/local
gp < example1.gp
\end{verbatim}

\subsection{Detailed example for a Hilbert class field}

To illustrate the approach, we provide an example of a class polynomial.
We choose $K$ primitive such that $\Kr$ has odd class number
and $\Krzero$ has class number~$1$,
so that by \cite [Theorem~I.10.3]{Streng10} the constructed class field
is the Hilbert class field of~$\Kr$.

Let $K = \Q (x)$ be the primitive non-cyclic CM field with $x$ a root of
$X^4 + 57 X^2 + 661$ and let $\Oc = \OK$ be the maximal order of~$K$.
We have
$K_0 = \Q (y) \cong \Q (\sqrt 5)$ with $y = x^2$ a root of $Y^2 + 57 Y + 661$,
and $\OKzero = \Z [\omega]$ has narrow class number~$1$, where
$\omega = \frac {y + 34}{11}$ satisfies $\omega^2 - \omega - 1 = 0$.
A generator of the different is $\diffgen = 2 \omega - 1$,
which satisfies $\diffgen^2 = 5$.

We choose the CM type $\Phi = (\varphi_1, \varphi_2)$ as
\[
\varphi_1 (x) = i \, \sqrt {\frac {57 - 11 \sqrt 5}{2}},\quad
\varphi_2 (x) = i \, \sqrt {\frac {57 + 11 \sqrt 5}{2}},
\]
which implies
\[
\varphi_1 (\diffgen) = - \varphi_2 (\diffgen) = \sqrt 5,\quad
\varphi_1 (\omega) = \frac {1 + \sqrt 5}{2}\quad
\text { and }\quad
\varphi_2 (\omega) = \frac {1 - \sqrt 5}{2},
\]
where all square roots of real numbers are taken positive.

The reflex field of~$K$ is given by $\Kr = \Q (t)\subseteq \CC$ with
$t \approx 10.41248483930371 i$ a root of $X^4 + 114 X^2 + 605$;
it contains the real quadratic number
$\omegar = -\frac{1}{4} (t^{2} + 55) = \frac {1 + \sqrt {661}}{2},$
where the positive real square root has been taken.
We find that the Igusa class polynomial is
\begin {align*}
841 X^3 &+ (-5611098752 \omegar - 17741044214880) X^2 \\
        &+ (3232391784287232 \omegar - 68899837678801920) X \\
        &+ (7331944332391841792 \omegar - 131969791422849515520).
\end {align*}

The prime~$3$ is inert in $K_0$ and splits in $K / K_0$,
so by Theorem~\ref {th:n} we may choose $N=3$ and are assured of the
existence of a $z_1 \in K \backslash K_0$
representing a principally polarised abelian surface,
such that~$z_1$ is the root of a
quadratic polynomial $[A_1, B_1, C_1]$ (which we use from now on as a
short-hand notation for $A_1 X^2 + B_1 X + C_1$) over $\OKzero$
with $\gcd (A_1, 3) = 1$ and $3 \mid C_1$.
With the approach of Theorem~\ref{th:n}
we find $[A_1,B_1,C_1] = [1, 1, 3\omega + 6]$.

This quadratic form has discriminant
$D = -12 \omega - 23$ and a root
$z_1 = \frac {-x^3 - 34 x - 11}{22}$,
which can readily be verified to lead
to a $\xi_1$ as in Proposition~\ref {prop:twogenerators} that is positive
imaginary under the two embeddings $\varphi_1$ and $\varphi_2$.
Following Corollary~\ref {cor:symplecticbasisg2}, we obtain for $z_1$ the
period matrix
\[
\tau_1 \approx \left( \begin {array}{rr}
0.5 + 4.1498183124610 i & 0.5 + 1.8108031294328 i \\
0.5 + 1.8108031294328 i &       2.3390151830282 i
\end {array} \right)
.
\]
We compute
\[
f_1 = I_4 (\tau_1 / 3) / I_4 (\tau_1)
\approx
4.31041770567796242256320 -
1.05769871912283540433298 i,
\]
which is a class invariant by Theorem~\ref {th:maintheorem}.

The CM class group $\Cg_{\OK, \Phi} (1)$ is isomorphic to the image~$\Dg$
of the homomorphism
\begin{align}
\label{eq:completetypenorm}
\Cl(\OKr) &\rightarrow  \left\{ (\bg, \nu) :
\begin{array}{l}\bg \text { fractional ideal of } \OK,\\
\nu \in K_0, \nu \gg 0,
\Norm_{K/K_0} (\bg) = \nu \OKzero
\end{array}
\right\}
/ \sim,
\\
\ag &\mapsto (N_{\Phi^r,\OK}(\ag), N(\ag)),
\end{align}
where the equivalence relation $\sim$ is given by the subgroup
$\left\{(\mu \OK, \mu \cconj {\mu}) : \mu \in K^\times \right\}$;
it can be computed as in~\cite{EnTh14} by the function
\verb!shimura_group_type_norm_subgroup!
of~CMH.
In this particular example, the group $\Dg$ is of order 3.
We have implemented Algorithm~\ref {alg:nsystem} as the function
\verb!nsystem! in~CMH;
for $N = 3$, it outputs, besides the initial $[1, 1, 3 \omega + 6]$,
the polynomials
$$[7\omega + 24, 48\omega - 83, -72\omega + 117]
\text { and } [3\omega + 14, -12\omega + 49, -24\omega + 51].$$
As a quick check, one immediately sees that the $A$ are coprime to~$3$
and that the $B$ are congruent to~$1$ modulo~$6$; then it is readily
verified that the discriminants are the same.

The three period matrices are computed by the functions
\verb!symplectic_basis! and \verb!period_matrix! of~CMH.
They lead to the floating point polynomial
\begin{align*}
F(X) \approx &
\, X^3 + (-1520.81864577885822782322 + 358.629756234205144714067 i)X^2\\
& + (120340.426264405965468052 - 39203.2377567834013587592 i) X\\
& + (-454033.008835683648854405 + 276194.792435730065214643 i),
\end{align*}
the coefficients of which are (not necessarily integral) elements of the
reflex field~$\Kr$.
Using the command \verb!recognize_polynomial! of RECIP
(based on the LLL algorithm) on $\Kr$ and a more precise approximation
of~$F$ we find conjecturally
\begin {align*}
&d' F (X) =
2^3 \cdot 11^5 \cdot 31^2 \cdot X^3 \\
& + (8560748430 t^3 + 11670666480 t^2
+ 970800040530 t - 617685149664) X^2 \\
& + (401850769605 t^3 - 3039243175155 t^2
+ 38906895998175 t - 180513547604841) X \\
& + (- 2982488461975 t^3 + 4298737055525 t^2
- 290518295198065 t - 96097164139933)
\end {align*}
with $d' = 2^3\cdot 11^5 \cdot 31^2$.
Alternatively, we can also obtain this polynomial as follows.
Guessing the minimal polynomials of the coefficients
(using the GP command
\texttt {algdep ($\cdot$, 4)}, for instance) reveals a common denominator
of $d = 11^4 \cdot 31^2$; taking the index between the
polynomial order $\Z [t]$ and its integral closure $\OKr$ into account,
we multiply $F$ with $d' = 2^3 \cdot 11 \cdot d$.
Integral linear dependencies obtained by the GP command
\texttt {lindep} between each coefficient of $d' F$ and
$1$, $t$, $t^2$ and $t^3$ yield the class polynomial above.

The height of this polynomial is a bit smaller than that of the
classical polynomial obtained from Igusa invariants above, but only
moderately so. Moreover since the polynomial has coefficients in the
CM field~$\Kr$ instead of its totally real subfield~$\Krzero$, printing
all coefficients actually takes more space.
But maybe it is not very surprising that quotients of
Igusa invariants do not result in a substantial gain in size: They are an analogue in
dimension~$2$ of quotients of the elliptic modular form~$\Delta$, which
in turn are the $24$-th powers of $\eta$-quotients; only lower powers of such
quotients are known to yield smaller class invariants \cite {EnMo14}.

In our case, it turns out that the $\sqrt {f_i}$ also lie in the Hilbert
class field (and thus generate it). The ``reason'' for this is that $h_4$
is the square of a Hilbert modular form for $\OKzero$, a situation that
we will examine in a future article. The class polynomial with roots
$\sqrt {f_1}$, $\sqrt {f_2}$ and $-\sqrt {f_3}$ (where all square
roots are taken with positive real part) is
conjecturally given by
\begin {align*}
F =
& \; 2^3 \cdot 11^3 \cdot 31 \cdot X^3 \\
& + (   44850 t^3 -  26268 t^2 +  5007630 t -  13168716) X^2 \\
& + (- 639765 t^3 + 657855 t^2 - 68212395 t -  21782871) X \\
& + (  693935 t^3 - 453871 t^2 + 68999645 t + 182497403).
\end {align*}

\subsection {Real example with a ramified level}

The following example illustrates Theorem~\ref {th:realramified}
for getting class invariants with real class polynomials.
We will use the level $N=2$.
Let $K = \Q (x)$ be the non-Galois quartic CM field with $x$ a root of
$X^4 + 18 X^2 + 68$, which has real subfield $K_0 = \Q (\sqrt {13})$,
where $\sqrt{13} = x^2+9$.
Consider
again the maximal order $\Oc = \OK$.
The hypothesis \eqref{eq:cchyp} in Theorem~\ref{th:realramified}
holds by Proposition~\ref{prop:cchyp}\eqref{it:ccg2}.

The real subfield of the reflex field is
$\Krzero = \Q (\omegar)$ with $\omegar = \frac {1 + \sqrt {17}}{2}$,
with again the usual embedding taking a positive real square root.
Using CMH and RECIP again with the function~$i_1$,
we obtain the following Igusa class polynomial for~$K$:

{
\tiny
\begin {align*}
    &19^2 \cdot 59^2\cdot X^4
      + (1381745663216332313130 \omegar - 3547293859211493542130) X^3 \\
    & + (148473995403415029782069841975 \omegar - 380321923961391822525781469475) X^2 \\
    & + (5344671730358474048907677495421000 \omegar - 13690639163949002342762017668129000) X \\
    & + (52888480565700710835194263641602550000 \omegar - 135476567266153427225864462713788270000).
\end {align*}
}

The prime~$2$ is inert in $K_0 / \Q$ and ramified in $K / K_0$.
Let $\omega = \frac{x^2+10}{2} = \frac{1 + \sqrt{13}}{2}\in K$,
which generates~$K_0$.
We fix the initial form as
\[
A_1 = 1, \quad B_1 = 0, \quad C_1 = 16 \omega + 22,
\]
where $C_1$ is divisible by $N=2$. We will thus obtain
a class polynomial defined over $\Krzero$ by Theorem~\ref {th:realramified}.

The group $\Cg_{\OK, \Phi} (1)$ is cyclic of order~$4$, and
a $2$-system is computed as
\begin {align*}
& A_2 = A_3 = 5 \omega + 7,
\quad &B_2 &= -B_3 = 8 \omega + 8,
\quad &C_2 &= C_3 = 2 \omega + 8, \\
& A_4 = -\omega + 5,
\quad &B_4 &= 0,
\quad &C_4 &= 6 \omega + 8,
\end {align*}
so that $f (\tau_1)$ and $f (\tau_4)$ are real and $f (\tau_2)$ and
$f (\tau_3)$ are complex conjugates whenever $f$ is a function for $\Gamma^0 (2)$
obtained as a quotient of two forms with rational $q$-expansions.
For $f = i_1=h_4h_6/h_{10}$, we get exactly the polynomial
above.
For the function $f = X^2 / Y$ of \S\ref {ssec:thetaproducts}, we obtain
numerically the following class polynomial:
{
\begin {align*}
X^4
& -19506.96702413769684992872390543869231117 X^3 \\
& +34104.71087980584199704592143935514042024 X^2 \\
& -31621.31544923554295971286232200559204341 X \\
& +17775.00158513035457655426023489497853378,
\end {align*}
}
which can be rewritten over $K_0$ as
{
\begin {align*}
19^4 \cdot 59^2 \cdot X^4
&+ (41960216624328 \omegar - 116332595812008) X^3 \\
&+ (-924565238142480 \omegar + 2383794199841616) X^2 \\
&+ (8404908240715776 \omegar - 21543961272975360) X \\
&+ (-10331028745814016 \omegar + 26471539326320640),
\end {align*}
}

\noindent which is noticeably smaller than the Igusa class polynomial.
As $Y$ is the square of a modular form,
the function $f$ is the square of a modular function.
Again it turns out that the zeroes of the previous polynomial
are squares of the zeroes of an even smaller polynomial, namely of:
{
	\begin {align*}
	19^2\cdot 59\cdot X^4 &+ (-2523732\omegar + 9426660)X^3 + (17576244\omegar - 46804644)X^2\\
	&+ (-15869952\omegar + 38596608)X - 49372416\omegar + 129309696.
\end{align*}
}

\subsection {Real example with a double Igusa quotient}

The following example illustrates Theorem~\ref {th:involutionS}.
Let $\Oc$ be the maximal order of the non-Galois quartic
CM field $K = \Q (x)$ with $x$ a root of $X^4 + 53 X^2 + 601$ and with real subfield
$K_0 \cong \Q (\sqrt 5)$. The class group of $\OK$ is cyclic of order~$5$
and isomorphic to $\Cg_{\OK, \Phi} (1)$. The real subfield of the reflex
field is $\Krzero = \Q (\sqrt {601})$, and we identify the algebraic integer
$\omegar = \frac {1 + \sqrt {601}}{2}$ with its positive real embedding.
With CMH and RECIP and~$i_1$ we obtain the following Igusa
class polynomial for~$K$:

{
\tiny
\begin {align*}
&2^{40} \cdot 13^4 \cdot X^5 \\
&+ (-6140585422220204445794304 \omegar - 322904904921695447307780096) X^4 \\
&+ (-96632884032276403274175741952 \omegar - 4131427744203466842763320885248) X^3 \\
&+ (-961856435411091691207536138780672 \omegar - 19922426752533168631849612073238528) X^2 \\
&+ (-2810878875032206947279703590350876416 \omegar - 32507451628887950858017880191429021184) X^{\phantom{1}} \\
&+ (-3949991728992949515358757855080152530801 \omegar - 59187968308773159157484805661633506074674)\phantom{X^0}
\end {align*}
}

We fix $N = 6$, the product of two primes that are inert in $K_0 / \Q$ and
split in $K / K_0$. The hypothesis \eqref{eq:cchyp2}
then follows from Proposition~\ref{prop:cchyp2}\point {2}.
Moreover, by Theorem~\ref {th:n}, there is a quadratic
polynomial $A_1 X^2 + B_1 X + C_1$ representing a polarised ideal class
with $6 \mid C_1$ and $\gcd (C_1/6, 6) = 1$;
for instance, $A_1 = 1$, $B_1 = \omega + 5$ and $C_1 = 6 \omega + 12$,
where $\omega = \frac{1}{9}(x^2+31)$.
Let $z_1$ be a root of this
polynomial, and choose the CM type in a compatible way; finally let $\tau_1$
be the associated period matrix as in Corollary~\ref {cor:symplecticbasisg2}.
We consider the double Igusa quotient
\[
f = \frac {h_{10} (\tau / 2) h_{10} (\tau / 3)}{h_{10} (\tau) h_{10} (\tau / 6)}.
\]
Then by Theorem~\ref {th:maintheorem}, the value $f (\tau_1)$ is a class
invariant, and by Theorem~\ref {th:involutionS}, its minimal polynomial
is real.

A $6$-system containing this initial form $Q_1 = [A_1, B_1, C_1]$ is
computed as
\begin {align*}
A_1 &=                1, & B_1 &=       \omega +    5, & C_1 &=     6 \omega +    12; \\
A_2 &=  -3 \omega +   7, & B_2 &=   145 \omega -  211, & C_2 &= -1110 \omega +  1866; \\
A_3 &=   3 \omega +   7, & B_3 &= - 143 \omega +    5, & C_3 &=   294 \omega +   606; \\
A_4 &= -34 \omega + 157, & B_4 &= -3959 \omega + 2309, & C_4 &=  4194 \omega + 34356; \\
A_5 &=     \omega +   2, & B_5 &= -  71 \omega +    5, & C_5 &=   180 \omega +   546.
\end {align*}

Letting $\tau_i$ denote the associated period matrices obtained by
Corollary~\ref {cor:symplecticbasisg2}, the conjugate $f (\tau_2)$ of the
class invariant is real, while $f (\tau_1)$ and $f (\tau_4)$ on one hand and
$f (\tau_3)$ and $f (\tau_5)$ on the other hand are complex conjugate pairs.

Numerically the class polynomial is approximated by
{
\begin {align*}
X^5
&-277.27759072275568417 X^4
 +3131337.2766719955916 X^3 \\
&-6196803.8120055534180 X^2
-2658.4275679312005124 X
-1.
\end {align*}
}Multiplying by the guessed denominator and identifying the resulting
coefficients as elements of $\Z [\omega_r]$ leads to the class
polynomial
{
\begin {align*}
2^4 \cdot 13^4 \cdot X^5
&+ (-53182948 \omegar + 551780268) X^4 \\
&+ (22828729975 \omegar + 1139705021035) X^3 \\
&+ (-46035175179 \omegar - 2244489935231) X^2 \\
&+ (10035944 \omegar - 1342872664) X
- 2^4\cdot 13^4,
\end {align*}
}which is already considerably smaller than the classical Igusa polynomial.

As $h_{10}$ is the square of a modular form,
the function $f$ is the square of a modular function.
It turns out that the zeroes of the previous polynomial
are squares of the zeroes of the following even smaller polynomial:
{
	\begin {align*}
2^2\cdot 13^2 X^5 &+ (1326\omegar + 23894)X^4 + (8833\omegar + 1025477)X^3 \\ &+ (-14003\omegar - 1482307)X^2 + (-2040\omegar - 6080)X - 2^2\cdot 13^2.
\end{align*}
}

\subsection{Necessity of the second complex conjugation hypothesis}
\label{ssec:cchyp2}

We now give an example that does not satisfy~\eqref{eq:cchyp2},
so Theorem~\ref{th:involutionS} does not apply, and for which
the associated class polynomial is defined over~$\Kr$,
but not $\Krzero$.

Let $\Oc=\OK$ be the maximal order of the non-Galois quartic CM field
$K = \Q(x)$ with $x$ a root of $ x^4 + 11x^2 + 29$
and with real subfield $K_0 \cong \Q(\sqrt{5})$.
The class group of $K$ is of order $2$ and generated by
the ideal $\pg_2 = (2,  x^2 + x + 5)$.
For every CM type of $K$,
the abelian varieties with CM by $\OK$ of that type are, up to
$\overline{\Q}$-isomorphism and the action of $\Gal(K_0^r/\Q)$, the Jacobians
of two curves over $K_0^r$ given in \cite[Table~2B]{BouyerStreng}.
In particular, the CM class group $\Cg_{\OK,\Phi}(1)$ is trivial.

We fix $N=14$, the product of the two primes $2$ and $7$ that
are inert in $K_0/\Q$. The prime $2$ ramifies as $(2) = \pg_2^2$
with $\pg_2$ as above,
and $7$ splits as $(7) = \pg_7 \overline{\pg_7}$ with
the principal ideal $\pg_7 = (x^2 + x + 6)$.

Let $\omega = x^2+6$ and 
consider the quadratic polynomial $AX^2 + BX + C$
with $A=1$, $B = 12\omega - 2$, $C = 14(2\omega + 3)$.
Let $z$ be a root of this polynomial and choose the CM type~$\Phi$
in a compatible way.
From $z$ we compute $\bg = (1)$, $\xi = (2x^3+16x)^{-1}$,
and $\NNg = N\OK+Az\OK = \pg_2\pg_7$, which is non-principal.
As $\Cg_{\OK,\Phi}(1)$ is trivial, we get that
$\mu N_{\Phi^r,\OK}(\cg)$ is principal in \eqref{eq:Nhyp}
and hence never equal to~$\NNg$.
So we get that \eqref{eq:Nhyp} does not hold,
and hence by Proposition~\ref{prop:cchyp}
neither does the second complex conjugation hypothesis~\eqref{eq:cchyp2}.

Consider the simple $I_4$ quotient
$$h(\tau) = I_4(\tau/N) / I_4(\tau),$$
which has rational $q$-expansion coefficients and is modular for $\Gamma^0(N)$.
The function $f = h + h\circ \iota$ is stable under the Fricke
involution $\iota : \tau \mapsto -N/\tau^{-1}$.

We use interval arithmetic and \cite{recip-3.4.1}
to evaluate $f$ in the point $\tau$ obtained from~$z$.
This yields
$$
f (\tau) \approx 
-65577.5 + 546773.0i, $$
with an error of absolute value less than $10^{-1}$,
which proves that $f(\tau)$
is non-real, and hence
confirms that its (linear) class polynomial does not have coefficients in~$K_0^r$.
This proves that condition \eqref{eq:cchyp2} is necessary in Theorem~\ref{th:involutionS}.

We used the complicated function $f$ instead of a double Igusa quotient
for the following reason.
Let $z_7 = z/7$, which is a root of $(7A)X^2 + BX + (C/7)$
and gives rise to $\tau_7 = \tau/7$.
Then the double Igusa quotient
$$\frac{h_k(\tau/2)h_k(\tau/7)}{h_k(\tau)h_k(\tau/14)}
= \frac{h_k(\tau/2)}{h_k(\tau)}\cdot
\left(\frac{h_k(\tau_7/2)}{h_k(\tau_7)}\right)^{-1}$$
is a quotient of simple Igusa quotients with ramified level~$2$,
which for this example both lie in~$K_0^r$ by Theorem~\ref{th:realramified},
hence do not illustrate the necessity of~\eqref{eq:cchyp2}.
In fact, as $\pg_7$ is principal, we find that both simple Igusa quotients
take the same value and hence this double Igusa quotient
is equal to~$1$.

\begin {sloppypar}
\printbibliography
\end {sloppypar}

\end {document}